\documentclass[letterpaper,11pt]{article}

\usepackage[T1]{fontenc}
\usepackage{lmodern}
\usepackage{amsmath,amssymb,mathtools,mathrsfs,bm}
\usepackage{amsthm}
\usepackage{aliascnt}
\usepackage[shortlabels]{enumitem}
\usepackage{booktabs,tabularx,array}
\usepackage[margin=1in]{geometry}
\usepackage{xcolor}
\definecolor{darkblue}{rgb}{0,0,.5}
\usepackage[breaklinks=true,pdfstartview=FitH]{hyperref}
\hypersetup{colorlinks=true,linkcolor=darkblue,citecolor=darkblue,urlcolor=darkblue}
\usepackage[capitalize,nameinlink]{cleveref}
\usepackage{microtype}
\usepackage[numbers,sort&compress]{natbib}

\allowdisplaybreaks
\AtBeginDocument{%
  \setlength{\abovedisplayskip}{5pt}
  \setlength{\belowdisplayskip}{5pt}
  \setlength{\abovedisplayshortskip}{3pt}
  \setlength{\belowdisplayshortskip}{3pt}
  \setlength{\jot}{1pt}
}

\newtheorem{theorem}{Theorem}[section]
\newaliascnt{lemma}{theorem}
\newtheorem{lemma}[lemma]{Lemma}
\aliascntresetthe{lemma}
\newaliascnt{proposition}{theorem}
\newtheorem{proposition}[proposition]{Proposition}
\aliascntresetthe{proposition}
\newaliascnt{corollary}{theorem}
\newtheorem{corollary}[corollary]{Corollary}
\aliascntresetthe{corollary}
\newaliascnt{fact}{theorem}

\aliascntresetthe{fact}
\newaliascnt{assumption}{theorem}

\aliascntresetthe{assumption}
\theoremstyle{definition}
\newaliascnt{definition}{theorem}
\newtheorem{definition}[definition]{Definition}
\aliascntresetthe{definition}
\newaliascnt{example}{theorem}
\newtheorem{example}[example]{Example}
\aliascntresetthe{example}
\newaliascnt{remark}{theorem}
\newtheorem{remark}[remark]{Remark}
\aliascntresetthe{remark}

\crefname{assumption}{Assumption}{Assumptions}
\crefname{theorem}{Theorem}{Theorems}
\crefname{lemma}{Lemma}{Lemmas}
\crefname{proposition}{Proposition}{Propositions}
\crefname{corollary}{Corollary}{Corollaries}
\crefname{fact}{Fact}{Facts}
\crefname{definition}{Definition}{Definitions}
\crefname{example}{Example}{Examples}
\crefname{remark}{Remark}{Remarks}

\newcommand{\Hh}{\mathcal H}

\newcommand{\Rr}{\mathcal R}

\newcommand{\Id}{\operatorname{Id}}
\newcommand{\ran}{\operatorname{ran}}
\newcommand{\dom}{\operatorname{dom}}
\newcommand{\gra}{\operatorname{gra}}
\newcommand{\zer}{\operatorname{zer}}
\newcommand{\Fix}{\operatorname{Fix}}
\newcommand{\dist}{\operatorname{dist}}
\newcommand{\cl}{\operatorname{cl}}

\newcommand{\N}{\mathbb N}
\newcommand{\R}{\mathbb R}
\newcommand{\ip}[2]{\left\langle #1,#2\right\rangle}
\newcommand{\norm}[1]{\left\lVert #1\right\rVert}

\newcommand{\weakto}{\rightharpoonup}

\newcommand{\eps}{\varepsilon}

\title{Degenerate Preconditioned Proximal Point Algorithms as Forward Methods for Nonmonotone Inclusions}
\author{Jan Harold Alcantara\thanks{\url{janharold.alcantara@riken.jp}.  Center for Advanced Intelligence Project, RIKEN, Tokyo, Japan.}}
\date{\today}

\begin{document}
\maketitle

\begin{abstract}
	We study degenerate preconditioned proximal point algorithms
	(PPPA) for finding zeros of nonmonotone set-valued operators in Hilbert
	spaces under solution-relative $(U,V)$-semimonotonicity, with
	possibly indefinite $U$ and $V$. With a
	closed-range preconditioner $P$, the reduced PPPA dynamics
	admit a forward representation through the complementary
	resolvent of a parallel composition, without requiring
	monotonicity or admissibility. Relative semimonotonicity of
	the range-restricted operator transfers to the reduced
	forward map, allowing our forward theory for possibly
	set-valued mappings to yield Fej\'er estimates and residual
	rates. 	Under suitable conditions, including demiclosedness of the
	reduced map at zero, the projected PPPA iterates and shadows
	converge weakly to a common projected zero. Without
	demiclosedness, a counterexample shows that the projected
	sequences may converge to a point outside the projected zero
	set, even when the original operator is continuous. Bounded
	recovery condition yields weak convergence of the full (non-projected) shadows.
	The monotonicity modulus $U$ yields additional results:
	under suitable conditions, $U\succeq\tau P$ with $\tau>0$
	gives strong and linear convergence of the projected PPPA
	sequences, while $U\succeq-\sigma P$ with $0\le\sigma<1/2$
	still permits local linear rates under metric subregularity.
These strong convergence and linear rate guarantees
extend to the full shadows
when a Lipschitz recovery estimate holds, such as when $U$ is strongly monotone on $\ker P$. For Chambolle--Pock and Douglas--Rachford
	splitting, we establish Lipschitz recovery even when
	$U=0$ by exploiting the primal--dual structure.
	Under suitable parameter conditions, full-shadow
	convergence guarantees, including $R$-linear rates,
	extend to the original PPPA iterates.
\end{abstract}



\section{Introduction}
Let $\Hh$ be a real Hilbert space and let
$A:\Hh\rightrightarrows\Hh$. We study the relaxed preconditioned
proximal point algorithm (PPPA)
\begin{equation}
	\bar x^k\in(P+A)^{-1}Px^k,\qquad
	x^{k+1}=x^k+\lambda_k(\bar x^k-x^k),
	\label{eq:intro-pppa}
\end{equation}
for finding a zero of $A$, where $P$ is bounded, self-adjoint, and
positive semidefinite. We call $\{\bar x^k\}$ the \emph{shadow sequence}.
For $P=\gamma^{-1}\Id$ and $\lambda_k\equiv1$, this is the classical
proximal point method of Martinet \cite{Martinet1970} and Rockafellar
\cite{Rockafellar1976}. Degenerate preconditioners arise
naturally in operator splitting: Douglas--Rachford splitting,
Chambolle--Pock at the critical stepsize, and graph-based splitting
schemes \cite{BrediesChencheneNaldi2024} admit such representations. Further examples of such representations are given in \cite{EvensEtAl2025}.

A systematic Hilbert-space theory of degenerate PPPA for monotone inclusions was developed by
Bredies--Chenchene--Lorenz--Naldi \cite{BrediesEtAl2022} through a
reduction to the range of $P$ and a parallel composition. Their framework also accommodates monotonicity restricted to graph points whose operator values lie in $\ran P$, a condition further studied by Xue--Zhang \cite{XueZhang2026}. Beyond monotonicity, nondegenerate proximal methods
have been studied under cohypomonotonicity
\cite{CombettesPennanen2004}, negative and
star-negative comonotonicity \cite{GorbunovEtAl2023},
and submonotonicity \cite{LukeTam2025}.
A local nonmonotone PPPA analysis with matrix-valued
relaxation also yields extensions of Spingarn's method
and progressive decoupling
\cite{EvensLatafatPatrinosSpingarn}.
For degenerate PPPA,
Evens--Pas--Latafat--Patrinos \cite{EvensEtAl2025}
developed a finite-dimensional analysis under an
oblique weak Minty condition, with applications to
Douglas--Rachford and Chambolle--Pock
\cite{EvensLatafatPatrinos2025}.

A common basis for comparing these assumptions is the
framework of $(U,V)$-semimonotonicity
\cite{AlcantaraDaoTakeda2026,EvensEtAl2025,EvensLatafatPatrinos2025}.
Here, $U$ and $V$ are self-adjoint operators, which we refer to as the
\emph{monotonicity modulus} and \emph{comonotonicity modulus},
respectively. The pairwise condition requires
$\langle x-y,a-b\rangle\ge
\langle x-y,U(x-y)\rangle+\langle a-b,V(a-b)\rangle$
for every $(x,a),(y,b)\in\gra A$. In this work, we consider its solution-relative counterpart, which requires the
inequality only when $(y,b)=(x^\star,0)$ for prescribed
$x^\star$ in a nonempty set $S^\star\subseteq\zer A$. 
The solution-relative condition with $U=0$ recovers
the $V$-oblique weak Minty condition studied in the
finite-dimensional setting in \cite{EvensEtAl2025},
which need not imply monotonicity even when $V=0$.
Allowing a negative monotonicity modulus $U$ further
enlarges this class.
\Cref{tab:literature-comparison} summarizes the operator
assumptions and convergence guarantees of representative
PPA and PPPA analyses.

\begin{table}[htbp]
	\centering
	\caption{Comparison of operator assumptions and convergence
		guarantees for PPA and PPPA.}
	\label{tab:literature-comparison}
	\begingroup
	\small
	\setlength{\tabcolsep}{2pt}
	\renewcommand{\arraystretch}{1.16}
	\begin{tabular}{@{}
			>{\raggedright\arraybackslash}p{.192\linewidth}
			>{\centering\arraybackslash}p{.175\linewidth}
			>{\centering\arraybackslash}p{.135\linewidth}
			>{\raggedright\arraybackslash}p{.088\linewidth}
			>{\raggedright\arraybackslash}p{.078\linewidth}
			>{\raggedright\arraybackslash}p{.100\linewidth}
			>{\raggedright\arraybackslash}p{.077\linewidth}
			>{\raggedright\arraybackslash}p{.077\linewidth}@{}}
		\toprule
		& \multicolumn{4}{c}{Assumptions and setting}
		& \multicolumn{3}{c}{Convergence guarantees} \\
		\cmidrule(lr){2-5}\cmidrule(l){6-8}
		Reference & $U$ & $V$
		& \shortstack[l]{Pairwise\\req.?}
		& \shortstack[l]{Degen.\\$P$?}
		& \shortstack[l]{Proj.\\conv.}
		& \shortstack[l]{Full\\conv.}
		& \shortstack[l]{Linear\\conv.} \\
		\midrule
		Rockafellar (1976) 
		& $0$ & $0$
		& $\checkmark$ & $\times$
		& $\checkmark$ & $\checkmark$ & $\checkmark$ \\
		Bredies et al.\ (2022)
		& $\alpha P,\ \alpha\ge0$ & $0$
		& $\checkmark^{*}$ & $\checkmark$
		& $\checkmark$ & $\checkmark$ & $\checkmark^{\dagger}$ \\
		Evens et al.\ (2025)$^{\ddagger}$
		& $0$ & general
		& $\times$ & $\checkmark$
		& $\checkmark^{\S}$ & $\checkmark$ & $\checkmark$ \\
		Luke--Tam (2025)$^{\ddagger,\P}$
		& $-\vartheta\Id$ & $-\vartheta\Id$
		& $\checkmark$ & $\times$
		& $\times$ & $\times$ & $\checkmark$ \\
		Xue--Zhang (2026)
		& $0$ & $0$
		& $\checkmark^{*}$ & $\checkmark$
		& $\checkmark$ & $\checkmark$ & $\times$ \\
		\textbf{Present work}
		& $U\succeq-\sigma P$ & general
		& $\times$ & $\checkmark$
		& $\checkmark$ & $\checkmark$ & $\checkmark$ \\
		\bottomrule
	\end{tabular}
	\endgroup
	
	\smallskip
	\begin{minipage}{\linewidth}
		\footnotesize\raggedright
		$^{*}$ The pairwise condition is restricted to graph
		points whose operator values lie in $\ran P$.\par
		$^{\dagger}$ The linear rate in
		\cite[Theorem~2.19]{BrediesEtAl2022}
		concerns the reduced algorithm.\par
		$^{\ddagger}$ Finite-dimensional setting.\par
		$^{\S}$ The projected solution-identification assertion
		of Evens et al. requires additional regularity;
		see \cref{ex:projected-counterexample}.\par
		$^{\P}$ For Luke--Tam, the displayed moduli refer to
		$\gamma A$. Their pairwise submonotonicity
		condition is local.
	\end{minipage}
\end{table}

\subsection{Approach and main ideas}
Building on the reduced-space formulation of
\cite{BrediesEtAl2022}, we study degenerate PPPA for
possibly nonmonotone inclusions by treating its reduced
dynamics as a forward method.
 When $P$ has closed range
$\Rr\coloneqq\ran P$, the factorization $P=CC^*$
gives reduced variables $z^k\coloneqq C^*x^k$
and residuals $q^k\coloneqq C^*(x^k-\bar x^k)$
satisfying
\begin{equation*}
	z^{k+1}=z^k-\lambda_kq^k,\quad \text{where} \quad q^k \in F(z^k),
\end{equation*}
and $F$ is the complementary resolvent of a
parallel composition of $A$ by $C^*$. This representation
holds as an identity of set-valued mappings,
\textit{without requiring monotonicity or admissibility} (see \cref{prop:parallel-composition}).
Moreover, relative semimonotonicity of either $A$ (or its
range-restricted graph
$A_{\Rr}(x)\coloneqq Ax\cap\Rr$) transfers to $F$ naturally (see \cref{lem:transfer-lemma}).
Together, the forward representation and this transfer property
allow the reduced dynamics to be analyzed through elementary
forward-iteration estimates. We develop these estimates for
possibly set-valued, relatively semimonotone mappings,
extending classical results for single-valued
quasi-cocoercive mappings (cf. \cref{def:UV-semi}(ii)) in
\cite{Maruster1977,DaoTamTruong2026}. It also establishes strong and linear convergence when the
monotonicity modulus of the forward map is strongly monotone, and
local linear rates under metric subregularity even when this
modulus is bounded below by a negative multiple of the identity.

The reduced variables $z^k=C^*x^k$ and
$z^k-q^k=C^*\bar x^k$ encode, respectively,
the projected PPPA iterates $\Pi_{\Rr}x^k$
and projected shadows $\Pi_{\Rr}\bar x^k$.
The forward theory yields Fej\'er monotonicity and,
under suitable parameter conditions, $q^k\to0$.
Consequently, both projected sequences are bounded,
and their difference converges strongly to zero.
However, these conclusions do not ensure that their
weak cluster points are projections of zeros of $A$, even when $A$ is continuous,
as shown by \cref{ex:projected-counterexample}.
This distinction motivates the study of demiclosedness
of $F$ at zero, which suffices for this projected solution
identification and, under an appropriate solution-set
condition, weak convergence of both projected sequences.
 We show that this condition on $F$
follows from demiclosedness of $A$ (or of the range-restricted operator $A_{\Rr}$) at zero
and the bounded-recovery condition (see \cref{def:recovery}(i)). These
conditions also ensure boundedness of $\{ \bar x^k\}$ and that their weak cluster points are
zeros of $A$. Weak convergence of $\{\bar x^k\}$ then follows
when $A$ has only one zero $x$ with $C^*x=z^\infty$, where
$z^\infty$ is the weak limit of $\{z^k\}$, as is the case under
admissibility (see \cref{def:admissible}). The conclusions extend
to the original PPPA iterates under an additional
relaxation condition. Thus, elementary forward estimates
settle the reduced dynamics, while projected
solution identification and full-shadow convergence require
separate demiclosedness and recovery arguments.

The monotonicity modulus $U$ provides further properties. Under
suitable conditions, the lower bound $U\succeq\tau P$ with $\tau>0$
translates into a strongly monotone modulus for $F$, so that the
general forward theory gives strong and linear convergence; a negative bound $U\succeq-\sigma P$ with $0\le\sigma<1/2$
similarly yields local linear rates under metric subregularity. Another role of $U$ is that its strong monotonicity on
$\ker P$ yields a Lipschitz recovery estimate (see
\cref{def:recovery}(ii)), which can be used to transfer strong convergence and
linear rates to the full PPPA shadows and, under suitable
relaxation, the iterates. For Chambolle--Pock and Douglas--Rachford
splitting, the primal--dual structure provides another route to
Lipschitz recovery, even when $U=0$.

\subsection{Contributions}
The main contributions are summarized as follows.

\begin{enumerate}[(i)]
\item \emph{Forward theory under relative semimonotonicity.}
We extend classical estimates for single-valued
quasi-cocoercive mappings to possibly set-valued,
solution-relative semimonotone forward maps.
Apart from Fej\'er estimates
and residual rates, we obtain strong convergence under a strongly monotone monotonicity modulus
(\cref{prop:forward-convergence}), and local linear rates
under metric subregularity when the monotonicity modulus has a negative lower bound
(\cref{prop:forward-linear}).

\item \emph{Projected convergence and solution identification.}
A two-dimensional counterexample with continuous $A$ shows
that the projected PPPA iterates and shadows can converge to a
point outside $\Pi_{\Rr}(\zer A)$
(\cref{ex:projected-counterexample}), demonstrating that the
projected identification and convergence assertions in
\cite[Theorem~2.4(iii),(v)]{EvensEtAl2025} require an
additional assumption. We establish projected solution
identification and weak convergence under demiclosedness of $F$
at zero (\cref{prop:projected-identification}), with sufficient
conditions involving demiclosedness of $A_{\Rr}$ and bounded
recovery (\cref{lem:sufficient-demiclosed}).
	
	\item \emph{Full-space recovery and convergence.}
Under bounded recovery, when $A$ has exactly one zero $x$ with
$C^*x=z^\infty$, where $z^\infty$ is the weak limit of $\{z^k\}$, we
establish weak convergence of the full PPPA shadows and, under
suitable relaxation, the iterates, without requiring admissibility
(\cref{thm:all-sequences}).
	Strong monotonicity of $U$ on $\ker P$, together with
	additional conditions, yields Lipschitz recovery
	(\cref{thm:kernel-lift}), transferring strong convergence
	and linear rates to the full shadows when
	$U\succeq\tau P$ with $\tau>0$
	(\cref{thm:weak-minty}).
	Local linear rates also extend to the full shadows when
	$U\succeq-\sigma P$ with $\sigma\in[0,1/2)$,
	under metric subregularity and suitable recovery conditions
	(\cref{cor:pppa-linear}).

\item \emph{Primal--dual splitting and recovery.}
For Chambolle--Pock and Douglas--Rachford splitting, full-graph
relative semimonotonicity with a block-diagonal $V$ yields Lipschitz
recovery even when $U=0$ (\cref{thm:pd-lift}). This shows that the
projected identification and full-shadow convergence results of
\cite{EvensEtAl2025,EvensLatafatPatrinos2025} hold, the latter
without continuity of the component resolvents (\cref{sec:cp}).
\end{enumerate}

\section{Preliminaries}\label{sec:prelim}

Throughout, $\Hh$ and $\mathcal K$ are real Hilbert spaces, and
$\N\coloneqq\{0,1,2,\ldots\}$. We denote by
$\mathcal B(\mathcal H,\mathcal K)$ the space of
bounded linear operators from $\mathcal H$ to $\mathcal K$, and
$\mathcal B(\mathcal H)\coloneqq\mathcal B(\mathcal H,\mathcal H)$.
We write $\Id_{\mathcal H}$ for the identity operator on
$\mathcal H$, and simply $\Id$ when the space is clear. A self-adjoint operator $U\in\mathcal B(\mathcal H)$
is said to be \emph{strongly monotone on a closed subspace}
$\mathcal S\subseteq\mathcal H$ if there exists $\mu>0$
such that $\ip{h}{Uh}\ge\mu\norm{h}^2$
for every $h\in\mathcal S$. For a closed subspace $\mathcal S$ of a Hilbert space,
$\Pi_{\mathcal S}$ denotes the orthogonal projector onto
$\mathcal S$. For a subset $S$ of a Hilbert space,
$\cl S$ denotes its closure and
$\dist(x,S)\coloneqq\inf_{s\in S}\norm{x-s}$. For $r>0$ and $\bar x\in \Hh$, we write
$B_r(\bar x)\coloneqq\{x\in\mathcal H:
\norm{x-\bar x}<r\}$. 
A sequence $\{z^k\}_{k\in\N}$ is said to be
\emph{Fej\'er monotone} with respect to a nonempty set $S$ if
$
\norm{z^{k+1}-z^\star}
\le\norm{z^k-z^\star}$
for all $z^\star\in S, ~k\in\N.$ Weak and strong convergence are denoted by $\weakto$ and $\to$,
respectively. For a set-valued operator $A$, we write
$\dom A\coloneqq\{x:Ax\ne\varnothing\}$,
$\ran A\coloneqq\bigcup_x Ax$,
$\gra A\coloneqq\{(x,a):a\in Ax\}$,
$\zer A\coloneqq\{x:0\in Ax\}$, and
$A^{-1}(a)\coloneqq\{x:a\in Ax\}$.
For a single-valued operator $T$, we write
$\Fix T\coloneqq\{x:Tx=x\}$.

\subsection{Semimonotonicity and closedness}
\label{subsec:prelim-semimonotone}

We define the class of semimonotone maps considered in recent works
\cite{AlcantaraDaoTakeda2026,EvensEtAl2025,EvensLatafatPatrinos2025},
together with its solution-relative counterpart.

\begin{definition}[$(U,V)$-semimonotonicity]
	\label{def:UV-semi}
	Let $A:\Hh\rightrightarrows\Hh$, and let
	$U,V\in\mathcal B(\Hh)$ be self-adjoint.
	\begin{enumerate}[(i)]
		\item We say that $A$ is
		\emph{$(U,V)$-semimonotone} if
		\begin{equation*}
			\ip{x-y}{a-b}
			\ge
			\ip{x-y}{U(x-y)}
			+\ip{a-b}{V(a-b)}
			\qquad
			\forall (x,a),(y,b)\in\gra A.
		\end{equation*}
		In particular, for $\beta>0$, we say that $A$
		is \emph{$\beta$-cocoercive} if it is
		$(0,\beta\Id)$-semimonotone.
		
		\item Let $\varnothing\ne S^\star\subseteq\zer A$.
		We say that $A$ is
		\emph{relatively $(U,V)$-semimonotone with respect
			to $S^\star$} if
		\[
		\ip{x-x^\star}{a}
		\ge
		\ip{x-x^\star}{U(x-x^\star)}
		+\ip{a}{Va}
		\qquad
		\forall (x,a)\in\gra A,\quad x^\star\in S^\star.
		\]
		When $U=0$, we also say that $A$ is
		\emph{weakly $V$-comonotone relative to $S^\star$},
		which is also known as the
		\emph{$V$-oblique weak Minty condition}
		\cite{EvensEtAl2025}.
		In particular, for $\beta>0$, we say that $A$
		is \emph{$\beta$-quasi-cocoercive relative to
			$S^\star$} if it is relatively
		$(0,\beta\Id)$-semimonotone with respect to
		$S^\star$.
	\end{enumerate}
\end{definition}

We refer to $U$ and $V$ as the
\emph{monotonicity modulus} and
\emph{comonotonicity modulus}, respectively,
of the semimonotonicity inequality. Global $(U,V)$-semimonotonicity implies relative
$(U,V)$-semimonotonicity with respect to every nonempty
subset of $\zer A$. We allow
$S^\star\subsetneq\zer A$ because the
solution-relative inequality may hold only at selected
zeros. 
\begin{remark}
	\label{rem:qco-interpretations}
	Quasi-cocoercivity admits several useful interpretations and
	consequences.
	\begin{enumerate}[(i)]
		
		\item Suppose that $A:\mathcal H\to\mathcal H$ is single-valued.
		If $S=\{x^\star\}$, then \cref{def:UV-semi}(ii) is precisely
		$1/\beta$-star-cocoercivity around $x^\star$ in the convention of
		\cite[Definition~2.1]{GorbunovLoizouGidel2022}. 	
		An equivalent condition appeared earlier as co-coercivity around
		a reference point in
		\cite[Definition~3.1]{LoizouEtAl2021}.
		If instead $S=\zer A$ and $T\coloneqq \Id-A$, then
		\cref{def:UV-semi}(ii) is precisely condition~(A) of
		\cite{Maruster1977}. It is also equivalent to conical
		$1/(2\beta)$-quasiaveragedness of $T$ in the sense of
		\cite[Definition~2.1]{DaoTamTruong2026}. 
		
		
		\item Unlike cocoercivity, quasi-cocoercivity does not force $A$
		to be single-valued, since it compares graph points only with the
		prescribed zeros. For example,
		$A:\Hh\rightrightarrows\Hh$ defined by $Ax\coloneqq \{x,2x\}$ is
		$1/2$-quasi-cocoercive relative to $\zer A=\{0\}$, although it is
		set-valued at every $x\ne0$. By contrast, cocoercivity forces
		single-valuedness on the domain. Nevertheless, quasi-cocoercivity constrains the values at $x^\star\in S$, since $a\in Ax^\star$ implies that
		$a=0$, i.e., $Ax^\star=\{0\}$. 
	\end{enumerate}
\end{remark}

We recall the notions of graph closedness and
demiclosedness at zero for set-valued operators;
see, e.g., \cite[Definition~2]{Isiogugu2013}
for the latter.

\begin{definition}[Graph closedness and demiclosedness at zero]	\label{def:closedness}
Let $A:\Hh\rightrightarrows\Hh$.
We say that $A$ is \emph{demiclosed at zero} if
$x^j\weakto x$, $a^j\to0$, and
$a^j\in Ax^j$ imply $0\in Ax$.
\end{definition}

%
%
%
%

\subsection{Degenerate metrics, admissibility, and recovery}
\label{subsec:prelim-degenerate}

Throughout this paper, $P\in\mathcal B(\Hh)\setminus\{ 0\}$ is a self-adjoint and positive semidefinite operator with closed range
$\Rr\coloneqq \ran P$. Then
$\Hh=\Rr\oplus\ker P$, and the restriction
$P|_{\Rr}:\Rr\to\Rr$ is bounded, self-adjoint, and strongly
positive. We equip $\Rr$ with the Hilbert structure inherited from $\Hh$
and use the canonical factorization
\begin{equation}
	C\coloneqq \iota_{\Rr}(P|_{\Rr})^{1/2}:\Rr\to\Hh,
	\qquad
	C^*=(P|_{\Rr})^{1/2}\Pi_{\Rr}:\Hh\to\Rr,
	\qquad
	CC^*=P,
	\label{eq:C-canonical}
\end{equation}
where $\iota_{\Rr}:\Rr\hookrightarrow\Hh$ is the canonical inclusion.
Thus, $C$ is injective, $\ran C=\Rr$, and $C^*$ is onto. Moreover,
$C^*|_{\Rr}=(P|_{\Rr})^{1/2}$ is a bounded isomorphism on $\Rr$, and hence
\begin{equation}
	\Pi_{\Rr}x=(P|_{\Rr})^{-1/2}C^*x,
	\qquad x\in\Hh.
	\label{eq:proj-via-C}
\end{equation}
We denote the seminorm
$\norm{x}_P^2\coloneqq \ip{x}{Px}=\norm{C^*x}^2$.

Following \cite[Definition~2.1]{BrediesEtAl2022}, we use the following
terminology.

\begin{definition}[Admissible preconditioner]
	\label{def:admissible}
	Let $A:\Hh\rightrightarrows\Hh$. A bounded self-adjoint
	positive semidefinite operator $P$ is an
	\emph{admissible preconditioner} for $A$ if
	$(P+A)^{-1}P:\Hh\to\Hh$ is single-valued
	and has full domain.
\end{definition}

Since $P=CC^*$, the preconditioned resolvent factors
through $\Rr$ as
\begin{equation}
	(P+A)^{-1}P=\Phi C^*,\quad \text{where} \quad 
	\Phi\coloneqq(P+A)^{-1}C:
	\Rr\rightrightarrows\Hh.
	\label{eq:shadow-recovery-map}
\end{equation}
We call $\Phi$ the \emph{shadow-recovery map}, since the
shadows in \eqref{eq:intro-pppa} satisfy
$\bar x^k\in\Phi(C^*x^k)$. As $C^*$ is onto, $P$
is admissible for $A$ if and only if $\Phi$ is
single-valued and has full domain. We also use two
recovery properties that do not require single-valuedness.

\begin{definition}[Bounded and Lipschitz recovery]
	\label{def:recovery} Let $\Phi$ be given by \eqref{eq:shadow-recovery-map}.
	\begin{enumerate}[(i)]
		\item The \emph{bounded-recovery condition}
		holds if $\Phi$ maps bounded sets of
		$\Rr$ to bounded sets of $\Hh$.
		
		\item Let $\varnothing\ne S^\star\subseteq
		\zer A$. A \emph{Lipschitz recovery estimate}
		holds relative to $S^\star$ if there exists
		$K>0$ such that
		\begin{equation}
			\norm{\bar x-x^\star}
			\le K\norm{z-C^*x^\star}
			\qquad
			\forall x^\star\in S^\star,\quad
			z\in\Rr,\quad \bar x\in\Phi(z).
			\label{eq:lipschitz-recovery}
		\end{equation}
		For $x^\star\in\zer A$, a
		\emph{local Lipschitz recovery estimate at
			$x^\star$} means that the same inequality in \eqref{eq:lipschitz-recovery} holds for every
		$z\in B_r(C^*x^\star)$ and
		$\bar x\in\Phi(z)$, for some $K,r>0$.
	\end{enumerate}
\end{definition}

The local Lipschitz recovery estimate is precisely
\textit{outer Lipschitz continuity} of $\Phi$ at
$C^*x^\star$, relative to $\dom\Phi$, in the
sense of \cite[Chapter 3D]{DontchevRockafellar2014}.
It does not require $\Phi$ to be single-valued
near $C^*x^\star$ or Lipschitz continuous
between arbitrary reduced points.

\begin{lemma}[Consequences of Lipschitz recovery]
	\label{lem:recovery-consequences}
	If \eqref{eq:lipschitz-recovery} holds relative
	to $\varnothing\ne S^\star\subseteq\zer A$,
	then the bounded-recovery condition holds and
	$\Phi(C^*x^\star)=\{x^\star\}$ for every
	$x^\star\in S^\star$. In particular, \(x^\star\) is the unique zero of \(A\) satisfying \(C^*x=C^*x^\star\).
The identity $\Phi(C^*x^\star)=\{x^\star\}$ and
uniqueness conclusion also hold under a local Lipschitz
recovery estimate at $x^\star$.
\end{lemma}

\begin{proof}
	Fix $x^\star\in S^\star$ and let $K$ be as in
	\eqref{eq:lipschitz-recovery}. For any bounded
	$B\subseteq\Rr$ and $\bar x\in\Phi(z)$ with
	$z\in B$,
	\[
	\norm{\bar x}\le\norm{x^\star}
	+K\bigl(\norm z+\norm{C^*x^\star}\bigr).
	\]
	It immediately follows that $\Phi(B)$ is bounded. Since
	$0\in Ax^\star$,
	$Px^\star=C(C^*x^\star)\in(P+A)x^\star$,
	so $x^\star\in\Phi(C^*x^\star)$.
	Applying \eqref{eq:lipschitz-recovery}
	at $z=C^*x^\star$ gives
	$\Phi(C^*x^\star)=\{x^\star\}$.
	Every $x\in\zer A$ with
	$C^*x=C^*x^\star$ belongs to this set
	and therefore equals $x^\star$. Both conclusions use the recovery estimate
	only at $z=C^*x^\star$, so they remain valid
	under local Lipschitz recovery.
\end{proof}

\section{Forward iterations under relative semimonotonicity}
\label{sec:forward-rel-semi}
Let $\mathcal{K}$ be a real Hilbert space. Given  a set-valued mapping
$F:\mathcal K\rightrightarrows\mathcal K$, we study forward iterations 
\begin{equation}
	z^{k+1} = z^k - \lambda_k q^k ,\quad \text{where} \quad q^k \in F(z^k),
	\label{eq:forward}
\end{equation}
under the assumption that $F$ is relatively
$(U,V)$-semimonotone with respect to a nonempty subset of its
zeros, where $V$ is strongly monotone. We first consider in \cref{sec:forward-convergence} the case $U\succeq0$,
which effectively reduces to quasi-cocoercive forward iterations (cf. \cref{def:UV-semi}(ii)), with some additional convergence results when $U$ is strongly monotone.
The case where $U$ is bounded below by a negative multiple
of the identity is treated in \cref{sec:forward-linear-convergence} under an additional error bound assumption.

\subsection{Convergence}\label{sec:forward-convergence}
The following proposition establishes the convergence of \eqref{eq:forward} under relative semimonotonicity. 

\begin{proposition}[Relative semimonotone forward iteration]
	\label{prop:forward-convergence}
	Let $F:\mathcal K\rightrightarrows\mathcal K$ be relatively
	$(U,V)$-semimonotone with respect to
	$\varnothing\ne S\subseteq\zer F$, where
	$U\succeq0$ and $V\succeq \eta \Id$ with $\eta>0$. Let $z^0\in\mathcal K$, and let $\{z^k\}$ be generated by \eqref{eq:forward} with
	$\lambda_k\in (0,2\beta_k]$, where $\beta_k>0$ satisfying  $\ip{q^k}{Vq^k}\ge\beta_k\norm{q^k}^2$.
	Set $\delta_k\coloneqq \lambda_k(2\beta_k-\lambda_k)$.
	Then the following hold.

	\begin{enumerate}[(i)]
		\item $\{z^k\}_{k\in\N}$ is Fej\'er monotone with respect to
		$S$.

		\item If $\liminf_{k\to\infty}\delta_k>0$, then
		$\sum_{k=0}^\infty\norm{q^k}^2<+\infty$ and
		$q^k\to0$.
		If, in addition, $F$ is demiclosed at zero,
		then every weak cluster point of $\{z^k\}$ belongs
		to $\zer F$. If moreover $S=\zer F$, then
		$z^k\weakto z^\star$ for some $z^\star\in\zer F$.


		\item For every $z^\star\in S$ and $N\in\N$,
		$
		\norm{
			\frac{\sum_{k=0}^N\lambda_kq^k}
			{\sum_{k=0}^N\lambda_k}
		}
		\le
		\frac{2\norm{z^0-z^\star}}
		{\sum_{k=0}^N\lambda_k}.
		$
		In particular, if $\liminf_{k\to\infty}\lambda_k>0$,
		the right-hand side is $O(N^{-1})$.
		Moreover, if $\delta\coloneqq \inf_k\delta_k>0$, then
		$
		\norm{\frac{1}{N+1}\sum_{k=0}^Nq^k}
		\le
		\frac{\norm{z^0-z^\star}}
		{\sqrt{\delta(N+1)}}
		=O(N^{-1/2}).
		$

		\item For every $N$ such that
		$\sum_{k=0}^N\delta_k>0$ and every $z^\star\in S$,
		$
		\min_{0\le k\le N}\norm{q^k}^2
		\le
		\frac{\norm{z^0-z^\star}^2}
		{\sum_{k=0}^N\delta_k}.
		$
		In particular, if $\delta_k\ge\kappa>0$ for all $k$,
		then the right-hand side is at most
		$\frac{\norm{z^0-z^\star}^2}{\kappa(N+1)}$.
		
			\item For every $z^\star\in S$,
		$
		\sum_{k=0}^\infty
		\lambda_k\norm{U^{1/2}(z^k-z^\star)}^2<+\infty.
		$
		If $\liminf_{k\to\infty}\lambda_k>0$, then
		$U^{1/2}(z^k-z^\star)\to0$.
		If $U\succeq\mu\Id$ for some $\mu>0$, then the
		weaker condition $\sum_{k=0}^\infty\lambda_k=+\infty$
		already ensures $z^k\to z^\star$ strongly.
		Moreover, if $\liminf_{k\to\infty}\lambda_k>0$,
		then the convergence is $Q$-linear.
	\end{enumerate}
\end{proposition}

\begin{proof}
	First, we note from the relative semimonotonicity assumption that 
	\begin{equation}
		\ip{q^k}{z^k-z^\star}
		\ge
		\ip{z^k-z^\star}{U(z^k-z^\star)}
		+\beta_k\norm{q^k}^2
		\qquad
		\forall z^\star\in S.
		\label{eq:relative-semimonotone-forward-lemma}
	\end{equation}
	From \eqref{eq:forward} and \eqref{eq:relative-semimonotone-forward-lemma}, we have for every $z^\star\in S$ 
	\begin{align}
		\norm{z^{k+1}-z^\star}^2
		\le{}&
		\norm{z^k-z^\star}^2
		-2\lambda_k\norm{U^{1/2}(z^k-z^\star)}^2
		-\lambda_k\ip{q^k}{(2V-\lambda_k\Id)q^k}
		\notag\\
		\le{}&
		\norm{z^k-z^\star}^2
		-2\lambda_k\norm{U^{1/2}(z^k-z^\star)}^2
		-\delta_k\norm{q^k}^2.
		\label{eq:fejermonotone-relative-semimonotone}
	\end{align}
	This proves (i). Summing \eqref{eq:fejermonotone-relative-semimonotone} from $k=0$ to $N$ yields
\begin{equation}
	2\sum_{k=0}^N
	\lambda_k\norm{U^{1/2}(z^k-z^\star)}^2
	+\sum_{k=0}^N\delta_k\norm{q^k}^2
	\le\norm{z^0-z^\star}^2,
	\label{eq:relative-forward-summability}
\end{equation}
which implies the first assertion in (ii).
Since $\{z^k\}$ is bounded by (i) and $q^k\to0$, demiclosedness of $F$ at zero implies that every weak cluster point belongs to $\zer F$.
If also $S=\zer F$, the distances
$\norm{z^k-z^\star}$ converge for every $z^\star\in\zer F$,
and Opial's argument gives the weak convergence of
$\{z^k\}$ to some point in $\zer F$. This complets the proof of (ii).

	For (iii), summing \eqref{eq:forward} gives
	$
	\sum_{k=0}^N\lambda_kq^k=z^0-z^{N+1}.
	$
	By Fej\'er monotonicity,
	\[
	\norm{z^0-z^{N+1}}
	\le
	\norm{z^0-z^\star}
	+\norm{z^{N+1}-z^\star}
	\le2\norm{z^0-z^\star},
	\]
	which proves the first estimate.
	If $\delta\coloneqq \inf_k\delta_k>0$, Jensen's inequality
	and \eqref{eq:relative-forward-summability} give
	$
	\norm{\frac{1}{N+1}\sum_{k=0}^Nq^k}^2
	\le
	\frac{1}{N+1}\sum_{k=0}^N\norm{q^k}^2
	\le
	\frac{\norm{z^0-z^\star}^2}{\delta(N+1)}.
	$
	This proves (iii). On the other hand, (iv) follows immediately since 
	$
	\left(\min_{0\le k\le N}\norm{q^k}^2\right)
	\sum_{k=0}^N\delta_k
	\le
	\sum_{k=0}^N\delta_k\norm{q^k}^2$.
	
The summability assertion  in (v) follows from \eqref{eq:relative-forward-summability}.
	If $\liminf_k\lambda_k>0$, the convergence
	$U^{1/2}(z^k-z^\star)\to0$ follows. 
	If $U\succeq\mu\Id$ for some $\mu>0$, then
	$\sum_k\lambda_k\norm{z^k-z^\star}^2<+\infty$.
	By (i), the sequence
	$a_k\coloneqq \norm{z^k-z^\star}^2$ is nonincreasing,
	so $a_k\to a$ for some $a\ge0$.
	If $a>0$, then
	$
	\sum_k\lambda_ka_k
	\ge a\sum_k\lambda_k=+\infty,
	$
	a contradiction. Hence $a=0$, proving
	$z^k\to z^\star$ strongly. 
	Furthermore, if $\liminf_k\lambda_k>0$, then
	\eqref{eq:fejermonotone-relative-semimonotone} and
	$U\succeq\mu\Id$ give
	$
	\norm{z^{k+1}-z^\star}^2
	\le
	(1-2\mu\lambda_k)\norm{z^k-z^\star}^2.
	$
	Since $\lambda_k$ is eventually bounded away from zero,
	this yields $Q$-linear convergence.
	This completes the proof of (v). 
\end{proof}

For $q^k\ne0$, the largest admissible coefficient in
\cref{prop:forward-convergence} is $
\beta_k
=
\ip{q^k}{Vq^k}/\norm{q^k}^2 \geq \eta .
$
Thus, the relaxation adapts to the residual direction and may
exceed $2\eta$. We may also take $\beta_k\equiv\eta$,
which gives the fixed stepsize condition $0<\lambda_k<2\eta$.

\begin{remark}[Relation with classical fixed-point iterations]
	\label{rem:DTT}
	Suppose that $F:\mathcal K\to\mathcal K$ is single-valued
	and $\beta$-quasi-cocoercive relative to $S=\zer F$,
	corresponding to $U=0$ and $V=\beta\Id$.
	Taking $\beta_k\equiv\beta$, $T\coloneqq \Id-F$, and
	$\rho\coloneqq 1/(2\beta)$, we have $\delta_k=\rho^{-1}\lambda_k(1-\rho\lambda_k)$ and the iteration in
	\cref{prop:forward-convergence} becomes
	$
	z^{k+1}=(1-\lambda_k)z^k+\lambda_kTz^k.
	$
	By \cref{rem:qco-interpretations}(i), $T$ is
	conically $\rho$-quasiaveraged.
	In this setting, \cref{prop:forward-convergence}(i),(ii)
	recover the Fej\'er monotonicity and residual convergence
	in \cite[Proposition~2.2(i),(ii)]{DaoTamTruong2026}.
	The ergodic estimate in
	\cite[Proposition~2.2(iii)]{DaoTamTruong2026}
	also follows from the residual summability in
	\cref{prop:forward-convergence}(ii) and Jensen's inequality.
	Moreover, under demiclosedness of $F$ at zero,
	the weak convergence conclusion in
	\cref{prop:forward-convergence}(ii) is consistent with
	the classical result of \cite[Theorem~1]{Maruster1977}.
	
	Unlike these single-valued formulations,
	\cref{prop:forward-convergence} permits arbitrary selections
	from a set-valued $F$, operator-valued moduli
	$U,V$, and iteration-dependent coefficients $\beta_k$.
The nonnegative $U$-term additionally yields further conclusions in
\cref{prop:forward-convergence}(v).

\end{remark}

\subsection{Linear convergence under metric subregularity}
\label{sec:forward-linear-convergence}

Metric subregularity is a standard local error-bound property for
set-valued mappings; see, e.g.,
\cite[Section~3H]{DontchevRockafellar2014}.
Together with suitable descent properties, metric
subregularity has been used to establish local linear
convergence of proximal and fixed-point iterations,
both in Hilbert-space monotone settings
\cite{Leventhal2009} and in finite-dimensional
nonmonotone settings
\cite{LukeThaoTam2018,EvensEtAl2025}. 

\begin{definition}[Metric subregularity]
	\label{def:metric-subreg}
	Let $F:\mathcal K\rightrightarrows\mathcal K$ and let
	$(\bar z,\bar q)\in\gra F$. We say that $F$ is
	\emph{metrically subregular at $\bar z$ for $\bar q$}
	if there exist $\chi,R>0$ such that
	\[
	\dist(z,F^{-1}(\bar q))
	\le\chi\dist(\bar q,F(z))
	\qquad\forall z\in  B_R(\bar{z}).
	\]
\end{definition}

For $\bar q=0$, this is the local error bound
$\dist(z,\zer F)\le\chi\dist(0,F(z))$.
Equivalently, metric subregularity of $F$ corresponds to
calmness of $F^{-1}$; see \cite{DontchevRockafellar2014}. We establish local linear convergence under metric subregularity
of $F$ or $G\coloneqq(F^{-1}-\Id)^{-1}$.
Note that $q\in F(z)$ if and only if $q\in G(z-q)$,
and $\zer G=\zer F$.

\begin{proposition}[Local linear convergence under metric subregularity]
	\label{prop:forward-linear}
	Let $F:\mathcal K\rightrightarrows\mathcal K$ be relatively
	$(U,V)$-semimonotone with respect to $\zer F\ne\varnothing$,
	where $U\succeq-\eps\Id$ for some $\eps\ge0$ and
	$V\succeq\eta\Id$ for some $\eta>0$.
	Fix $z^\star\in\zer F$ and
	suppose that, for some $R>0$, $F$ has nonempty values
	on $B_R(z^\star)$ and one of the following holds:
	\begin{enumerate}[(a)]
		\item $F$ is metrically subregular at $z^\star$ for zero
		with modulus $\chi>0$, with the error bound holding
		on $B_R(z^\star)$;
		\item $G$ is metrically subregular at $z^\star$
		for zero with modulus
		$\widehat\chi>0$, with the error bound holding
		on $B_R(z^\star)$.
	\end{enumerate}
	Let $\{z^k\}$ be generated by \eqref{eq:forward} with
	$0<\lambda_k\le2\beta_k$, where $\beta_k\ge\eta$ satisfies
	$\ip{q^k}{Vq^k}\ge\beta_k\norm{q^k}^2$.
	Under (a) or (b), respectively, set
	\begin{equation}
		\Delta_k\coloneqq
		\begin{cases}
			\lambda_k(2\beta_k-\lambda_k-2\eps\chi^2),
			& \text{under (a)},\\
			\lambda_k[2\beta_k-\lambda_k
			-2\eps(\widehat\chi+|1-\lambda_k|)^2],
			& \text{under (b)}.
		\end{cases}
		\label{eq:local-gap}
	\end{equation}
	and suppose that $\Delta_k\ge\kappa$ for all $k$, where $\kappa >0$. Define $m_\eps(t)\coloneqq(1+\sqrt{1+4t\eps})/(2t)$.
	Then there exist $r\in(0,R)$, $\rho\in[0,1)$, and
	$\Gamma>0$ (depending only on $R$, $\kappa$, $\eps$,
	$\eta$, $\norm V$, and $\chi$ or $\widehat\chi$) such that,
	whenever $z^0\in B_r(z^\star)$, the following hold.
	\begin{enumerate}[(i)]
		\item The iterates satisfy $z^k\in B_R(z^\star)$
		for all $k$, and, under (b), also
		$z^k-q^k\in B_R(z^\star)$. 
		\item If $d_k\coloneqq\dist(z^k,\zer F)$, then $d_{k+1}\le\rho d_k$, and hence
		$d_k\le d_0\rho^k$. Moreover,
		$\norm{q^k}\le m_\eps(\eta)d_0\rho^k$, and
		$\norm{z^{k+1}-z^k}\le(1-\rho)\Gamma d_k$.
		
		\item The sequence $\{z^k\}$ converges strongly
		to some $z^\infty\in\zer F$, and, for every $m\ge k$,
		$\norm{z^m-z^k}\le\Gamma d_k$ and
		$\norm{z^k-z^\infty}\le\Gamma d_k$.
		
		\item The auxiliary sequence
		$\widehat z^k\coloneqq z^k-q^k$ converges strongly
		to $z^\infty$, with
		$\norm{\widehat z^k-z^\infty}
		\le(\Gamma+m_\eps(\beta_k))d_k$.
	\end{enumerate}

If $\eps=0$, the bounds in (iii)--(iv) hold with
$\Gamma$ replaced by $\min\{2,\Gamma\}$,
and $m_0(\beta_k)=1/\beta_k$.
\end{proposition}

\begin{proof}
	For every $z\in\zer F$, relative semimonotonicity
	and $V\succeq\eta\Id$ imply $F(z)=\{0\}$.
	Thus, if $q^k=0$ for some $k$, the iteration
	is stationary thereafter. We may therefore
	assume that $q^k\ne0$ for all $k$.
	For $q^k\ne0$, note that $\eta\le\beta_k\le\bar\beta$ where $\bar\beta\coloneqq\norm V$.
	Relative semimonotonicity gives, for every
	$w\in\zer F$,
	\begin{equation}
		\ip{q^k}{z^k-w}
		\ge\beta_k\norm{q^k}^2
		-\eps\norm{z^k-w}^2.
		\label{eq:violated-qco}
	\end{equation}
	By Cauchy--Schwarz and solving the quadratic inequality for $\norm{q^k}$%
	\footnote{For $a,b\ge0$ and $t>0$,
		$ta^2\le ab+\eps b^2$ implies
		$a\le(1+\sqrt{1+4t\eps})b/(2t)$.},
	\eqref{eq:violated-qco} gives
	$\norm{q^k}\le m_\eps(\beta_k)\norm{z^k-w}$, and taking the infimum over $w\in\zer F$ yields
	$\norm{q^k}\le m_\eps(\beta_k)d_k$.
	Since $m_\eps$ is nonincreasing and
	$\eta\le\beta_k\le\bar\beta$, we have
	$\norm{q^k}\le Md_k$ where $M\coloneqq m_\eps(\eta)$. Moreover,
	\begin{equation}
		\norm{z^{k+1}-z^k}=\lambda_k\norm{q^k}
		\le2\beta_km_\eps(\beta_k)d_k
		=\bigl(1+\sqrt{1+4\eps\beta_k}\bigr)d_k
		\le\Theta d_k, \label{eq:diff-linear}
	\end{equation}
	where $\Theta\coloneqq1+\sqrt{1+4\eps\bar\beta}$.
	
	Set $\delta_k\coloneqq\lambda_k(2\beta_k-\lambda_k)$.
	From \eqref{eq:forward} and \eqref{eq:violated-qco},
	for every $w\in\zer F$,
	\begin{equation}
		\norm{z^{k+1}-w}^2
		\le(1+2\lambda_k\eps)\norm{z^k-w}^2
		-\delta_k\norm{q^k}^2.
		\label{eq:quasidescent}
	\end{equation}
	Taking the infimum over $w\in\zer F$ gives
	\begin{equation}
		d_{k+1}^2\le(1+2\lambda_k\eps)d_k^2
		-\delta_k\norm{q^k}^2.
		\label{eq:violated-distance}
	\end{equation}
	Under (a), if $z^k\in B_R(z^\star)$, the error bound
	gives $d_k\le\chi\norm{q^k}$, and hence
	$d_{k+1}^2\le(1-\Delta_k/\chi^2)d_k^2$.
	Under (b), if $z^k-q^k\in B_R(z^\star)$,
	then $q^k\in G(z^k-q^k)$ yields
	$\dist(z^k-q^k,\zer F)\le\widehat\chi\norm{q^k}$.
	Set $b_k\coloneqq\widehat\chi+|1-\lambda_k|$.
	The triangle inequality and \eqref{eq:forward} gives
	$d_{k+1}\le b_k\norm{q^k}$, so
	\eqref{eq:violated-distance} yields
	$(b_k^2+\delta_k)d_{k+1}^2
	\le(1+2\lambda_k\eps)b_k^2d_k^2$, and therefore $d_{k+1}^2 \leq (1-\Delta_k/(b_k^2+\delta_k))d_k^2$.
	As $\delta_k\le\bar\beta^2$,
	$b_k\le\widehat\chi+1+2\bar\beta$, and
	$\Delta_k\ge\kappa$,
	we have $d_{k+1}^2\leq  (1-\kappa/D)d_k^2$ where
	$D\coloneqq\bar\beta^2+
	(\widehat\chi+1+2\bar\beta)^2$.
	Thus, with
	\[\rho^2\coloneqq\begin{cases}
		\max\{0,1-\kappa/\chi^2\} , & \text{under (a)} \\
		\max\{0,1-\kappa/D\} ,& \text{under (b)} ,
	\end{cases} \]
	where $\rho\in[0,1)$,
	we have shown that 
	\begin{equation}
		d_{k+1}\le\rho d_k \text{ whenever } 
	\begin{cases}
		z^k\in B_R(z^\star), & \text{under (a), or} \\
		z^k-q^k\in B_R(z^\star) & \text{under (b)}.
	\end{cases}
	\label{eq:contraction-under-error-bound}
	\end{equation}

		Set $\Gamma\coloneqq\Theta/(1-\rho)$.
	Choose $r\in(0,R)$ so that $(1+\Gamma)r<R$
	under (a), or $(1+\Gamma+M)r<R$ under (b).
	Let $z^0\in B_r(z^\star)$. Under (a), we prove by induction that
$z^k\in B_R(z^\star)$ and $d_{k+1}\le\rho d_k$.
Indeed, suppose that $z^j\in B_R(z^\star)$
for all $j<k$. By
\eqref{eq:contraction-under-error-bound},
$d_j\le\rho^j d_0$ for $j\le k$.
Using \eqref{eq:diff-linear},
\begin{align}
	\norm{z^k-z^\star}
	&\le\norm{z^0-z^\star}
	+\Theta\sum_{j=1}^{k}d_{j-1}\notag \\
	&\le\norm{z^0-z^\star}+\Gamma d_0
	\le(1+\Gamma)\norm{z^0-z^\star}<R. \label{eq:invariance}
\end{align}
Thus, $z^k\in B_R(z^\star)$, and hence
$d_{k+1}\le\rho d_k$ by
\eqref{eq:contraction-under-error-bound}.
This completes the induction under (a).
		
Under (b), since $\norm{q^0}\le Md_0$, we have
	$\norm{z^0-q^0-z^\star}
	\le(1+M)\norm{z^0-z^\star}<R$.
	Thus, the base case $z^0-q^0 \in B_R(z^\star)$ holds.
	Suppose that $z^j-q^j\in B_R(z^\star)$
	for all $j<k$. By
	\eqref{eq:contraction-under-error-bound},
	$d_j\le\rho^j d_0$ for $j\le k$.
	Arguing as in (a), we obtain
	$\norm{z^k-z^\star}
	\le(1+\Gamma)\norm{z^0-z^\star}<R$.
	Moreover, since $\norm{q^k}\le Md_k\le Md_0$,
	\[
	\begin{aligned}
		\norm{z^k-q^k-z^\star}
		&\le\norm{z^k-z^\star}+\norm{q^k}\\
		&\le(1+\Gamma+M)\norm{z^0-z^\star}<R.
	\end{aligned}
	\]
	Thus, $z^k,z^k-q^k\in B_R(z^\star)$, and hence
	$d_{k+1}\le\rho d_k$ by
	\eqref{eq:contraction-under-error-bound}.
	This completes the induction under (b). Hence, we have proved (i) and (ii). 
	
	For every $m\ge k$, \eqref{eq:diff-linear} and triangle inequality gives
	$\norm{z^m-z^k}\le\Theta\sum_{j=k}^{m-1}d_j
	\le\Gamma d_k$.
	Since $d_k\to0$, the sequence $\{z^k\}$ is Cauchy
	and converges strongly to some $z^\infty\in\mathcal K$.
Letting $k\to \infty$ in \eqref{eq:invariance}, we have
	$\norm{z^\infty-z^\star}
	\le(1+\Gamma)\norm{z^0-z^\star}<R$.
	Moreover,
	$\dist(z^\infty,\zer F)
	\le\norm{z^\infty-z^k}+d_k\to0$,
	so $z^\infty\in\cl(\zer F)$.
	Choose $q\in F(z^\infty)$, which is possible
	since $z^\infty\in B_R(z^\star)$, and let
	$w_j\in\zer F$ converge to $z^\infty$.
	Relative semimonotonicity yields
	$\eta\norm q^2\le\ip{q}{z^\infty-w_j}
	+\eps\norm{z^\infty-w_j}^2$.
	Letting $j\to\infty$ gives $q=0$, and hence
	$z^\infty\in\zer F$.
	Letting $m\to\infty$ in 	$\norm{z^m-z^k}\le\Gamma d_k$ proves
	$\norm{z^k-z^\infty}\le\Gamma d_k$.
	This proves (iii).
	
	By (i), $q^k\to0$, and hence
	$\widehat z^k=z^k-q^k\to z^\infty$.
	The triangle inequality gives
	$\norm{\widehat z^k-z^\infty}
	\le[\Gamma+m_\eps(\beta_k)]d_k$,
	proving (iv).
	Finally, when $\eps=0$, \eqref{eq:quasidescent}
	shows that the sequence is Fej\'er monotone
	with respect to $\zer F$.
	Thus, for every $m\ge k$ and $w\in\zer F$,
	$\norm{z^m-z^k}\le
	\norm{z^m-w}+\norm{z^k-w}
	\le2\norm{z^k-w}$.
Taking the infimum over $w\in\zer F$ gives
$\norm{z^m-z^k}\le2d_k$.
Letting $m\to\infty$ yields the improved bounds
in (iii). The corresponding bound in (iv)
follows by the triangle inequality and
$m_0(\beta_k)=1/\beta_k$.
\end{proof}

Under (a), if $\eta>\eps\chi^2$, choosing
$\lambda_k\equiv\eta-\eps\chi^2$ and
$\beta_k\equiv\eta$ gives
$\Delta_k=(\eta-\eps\chi^2)^2>0$.
Under (b), $\lambda_k\equiv1$ and
$\beta_k\equiv\eta$ are admissible whenever
$2\eta-1>2\eps\widehat\chi^2$, since then
$\Delta_k=2\eta-1-2\eps\widehat\chi^2>0$.

Under the hypotheses of \cref{prop:forward-linear},
metric subregularity of $G$ at $z^\star$ for zero
implies that of $F$ with modulus $1+\widehat\chi$.
Indeed, for $z$ near $z^\star$ and $q\in F(z)$,
the residual bound gives
$\norm q\le m_\eps(\eta)\norm{z-z^\star}$,
so $z-q$ also lies in the corresponding neighborhood.
Hence, $\dist(z,\zer F)\le(1+\widehat\chi)\norm q$.
Thus, (b) implies (a) locally, although (b) can
yield a smaller contraction factor, as illustrated
in \cref{rem:ppa-comparison}.

If $U\succeq\mu\Id$ for some $\mu>0$,
then $\zer F=\{z^\star\}$, and relative
semimonotonicity and Cauchy--Schwarz give
$\mu\norm{z-z^\star}\le\norm q$
for every $(z,q)\in\gra F$.
Thus, $F$ is globally metrically subregular
with modulus $1/\mu$.
The corresponding direct convergence result is
given by \cref{prop:forward-convergence}(v).

\begin{remark}[Relation to almost averagedness and ordinary PPA]
	\label{rem:ppa-comparison}
	For constant $\lambda_k\equiv\lambda$ and
	$\beta_k\equiv\beta$, \eqref{eq:quasidescent}
	has the algebraic form of the pointwise
	almost-averaged inequality, with violation
	$2\lambda\eps$ and averaging constant
	$\lambda/(2\beta)$; see
	\cite[Proposition~2.1(iii)]{LukeThaoTam2018}.
	
	Let $t>0$ and consider the proximal point iteration
	$z^{k+1}\in J_{tA}z^k$, where
	$J_{tA}\coloneqq(\Id+tA)^{-1}$.
	This is \eqref{eq:forward} with $\lambda_k\equiv1$
	and $F\coloneqq\Id-J_{tA}$, since
	$q^k\coloneqq z^k-z^{k+1}\in F(z^k)$.
	Moreover, $G=(F^{-1}-\Id)^{-1}=tA$ and
	$\zer F=\zer A$. 
	Local submonotonicity of $tA$ with violation
	$\nu\ge0$ in the sense of
	\cite[Definition~2]{LukeTam2025} gives
	$\ip{q}{z-w}\ge\norm q^2-\nu\norm{z-w}^2$
	for nearby zeros $w$ and $q\in F(z)$
	satisfying the local graph restrictions. Thus, $F$ satisfies a local
	form of \eqref{eq:violated-qco} with $\beta_k\equiv 1$
	and $\eps=\nu$.
In the notation of \cref{def:UV-semi},
their condition is local
$(-\vartheta/t\,\Id,-t\vartheta\Id)$-semimonotonicity
of $A$, where
$\vartheta\coloneqq\nu/(1+2\nu)$;
cf.\ \cref{lem:transfer-lemma} with
$P=t^{-1}\Id$ and $\tau=-\vartheta$.
	
	Suppose that $A$ is metrically subregular at
	the reference zero for zero with modulus
	$\chi_A>0$, and set $a\coloneqq\chi_A/t$.
	Then $G=tA$ is metrically subregular with
	modulus $a$. For $\lambda_k=\beta_k=1$,
	the one-step estimate
	$d_{k+1}^2\le
	(1-\Delta_k/(b_k^2+\delta_k))d_k^2$
	in the proof of \cref{prop:forward-linear}(b)
	has $\delta_k=1$, $b_k=a$, and
	$\Delta_k=1-2\nu a^2$.
	Hence, the squared contraction factor is
	$\rho^2=(1+2\nu)a^2/(1+a^2)$, provided
	$2\nu a^2<1$.
	This is strictly smaller than the limiting
	squared factor $1+2\nu-(1+a)^{-2}$ in
	\cite[Theorem~2]{LukeTam2025}, whose
	contractivity condition $2\nu(1+a)^2<1$
	implies ours.
	When $\nu=0$, taking square roots gives
	the classical rate $\chi_A/\sqrt{\chi_A^2+t^2}$
	of \cite[Theorem~3.1]{Leventhal2009}. 
	
	This comparison concerns contraction factors rather than
convergence theorems. While \cref{prop:forward-linear}
requires a solution-relative inequality for all values
of $F$, \cite{LukeTam2025} assumes pairwise
submonotonicity on a restricted graph region, with
additional conditions ensuring the existence of
proximal steps within that region. Thus, neither
theorem directly subsumes the other.
\end{remark}

\section{Degenerate PPPA as a forward method}
\label{sec:pppa-forward}

We now specialize the forward theory of \cref{sec:forward-rel-semi} to
the PPPA
\begin{equation}
	\bar x^k\in(P+A)^{-1}Px^k,
	\qquad
	x^{k+1}=x^k+\lambda_k(\bar x^k-x^k).
	\label{eq:pppa}
\end{equation}
Throughout this
section, along with the standing assumptions in \cref{subsec:prelim-degenerate}, we assume that
$\dom((P+A)^{-1}P)=\Hh$, so that the shadow sequence $\{ \bar x^k\}$  is well defined.

\subsection{Reduced forward representation}
\label{subsec:pppa-forward}
We consider a reduced representation of
\eqref{eq:pppa} on $\Rr$ by setting
$z^k\coloneqq C^*x^k\in\Rr$, which we call the \textit{reduced PPPA sequence}. Observe that
$Px^k=Cz^k$, and so \eqref{eq:pppa} gives
\begin{equation}
	\bar x^k\in(P+A)^{-1}Cz^k,
	\qquad
	z^{k+1}
	=z^k-\lambda_k(z^k-C^*\bar x^k).\label{eq:pppa-reduced}
\end{equation}
This motivates the definition of the \textit{reduced forward map}
$F:\Rr\rightrightarrows\Rr$ by
\begin{equation}
	F(z)\coloneqq \{z-C^*\bar{x}:\bar{x}\in(P+A)^{-1}Cz\}.
	\label{eq:reduced-forward-map}
\end{equation}
Equivalently, $q\in F(z)$ if and only if
there exists $\bar{x}\in\Hh$ such that
$C^*\bar{x}=z-q$ and $Cq\in A\bar{x}$.
If $P$ is admissible for $A$, then
$F$ is single-valued.

\begin{proposition}[Reduced forward representation]
	\label{prop:forward-representation}
	Let $(x^k,\bar x^k)$ be generated by \eqref{eq:pppa},
	and set
	\begin{equation}
			z^k\coloneqq C^*x^k \quad \text{and} \quad q^k\coloneqq C^*(x^k-\bar x^k) = z^k-C^*\bar x^k.
		\label{eq:zq}
	\end{equation}
	Then
		\begin{align}
			Cq^k& =P(x^k-\bar x^k)\in A\bar x^k,
			\label{eq:graph-rel}\\
			z^{k+1}& =z^k-\lambda_kq^k,
			\quad \text{where} \quad 
			q^k\in F(z^k),\quad \text{and} 
			\label{eq:forward-main}\\
			\norm{q^k}& =\norm{x^k-\bar x^k}_P.
			\label{eq:residual-equality}
		\end{align}
\end{proposition}

\begin{proof}
	Since $\bar x^k\in(P+A)^{-1}Px^k$, we have
	$P(x^k-\bar x^k)\in A\bar x^k$. Then
	$Cq^k=CC^*(x^k-\bar x^k) = P(x^k-\bar x^k)\in A\bar x^k$. This proves \eqref{eq:graph-rel}. It also follows from the definition of $F$ that 
	$q^k=z^k-C^*\bar x^k\in F(z^k)$, which together with \eqref{eq:pppa-reduced}
	gives \eqref{eq:forward-main}.
	Finally, \eqref{eq:residual-equality} directly follows from \eqref{eq:zq}.
\end{proof}

The reduced forward map $F$ can also be expressed
in terms of the \emph{parallel composition of $A$ by $C^*$} \cite[Definition~25.39]{BauschkeCombettes2017}. For an operator $B$, we write
$J_B\coloneqq (\Id+B)^{-1}$ for its \emph{resolvent}, which need
not be single-valued, and denote by
\[
A_C\coloneqq C^*\triangleright A
\coloneqq (C^*A^{-1}C)^{-1}:\Rr\rightrightarrows\Rr
\]
the {parallel composition}.
Under appropriate monotonicity and admissibility
assumptions, the reduced-resolvent identity
$J_{A_C}=C^*(P+A)^{-1}C$ was established in
\cite[Theorem~2.13]{BrediesEtAl2022}.
The following proposition establishes these
identities without requiring monotonicity
or admissibility.

\begin{proposition}[Parallel-composition representation]
	\label{prop:parallel-composition}
	As identities of set-valued mappings on $\Rr$,
	\begin{align}
		& C^*(P+A)^{-1}C=J_{A_C}, \label{eq:parallel-resolvent}\\
		&F=J_{A_C^{-1}}=\Id_{\Rr}-J_{A_C},
		\label{eq:F-as-resolvent}\\
		&\zer F=\zer A_C=C^*(\zer A).
		\label{eq:zeros-F}
	\end{align}
\end{proposition}

\begin{proof}
	Let $z,p\in\Rr$ and $\bar{x}\in\Hh$ with $p=C^*\bar{x}$.
	Since $P\bar{x}=Cp$, we have
	$\bar{x}\in(P+A)^{-1}Cz$ if and only if
	$C(z-p)\in A\bar{x}$. Hence
	$p\in C^*(P+A)^{-1}Cz$ if and only if
	$p\in C^*A^{-1}C(z-p)$, that is,
	$z-p\in A_Cp$, or equivalently
	$p\in J_{A_C}(z)$. This proves \eqref{eq:parallel-resolvent}. On the other hand, by the definition of $F$, $q\in F(z)$ if and only if
	$z-q\in C^*(P+A)^{-1}Cz=J_{A_C}(z)$.
	Moreover, $z-q\in J_{A_C}(z)$ if and only if
	$q\in A_C(z-q)$, that is,
	$z-q\in A_C^{-1}q$, or equivalently
	$q\in J_{A_C^{-1}}(z)$.
	This proves \eqref{eq:F-as-resolvent}. Finally, $0\in F(z)$ if and only if
	$z\in J_{A_C}(z)$, that is,
	$0\in A_Cz$, or equivalently
	$z\in C^*A^{-1}(0)=C^*(\zer A)$.
\end{proof}

Let $p^k\coloneqq C^*\bar x^k=z^k-q^k $.
By \cref{prop:parallel-composition},
$p^k \in J_{A_C}(z^k)$, and hence
\eqref{eq:forward-main} can also be written as
$
z^{k+1}=(1-\lambda_k)z^k+\lambda_kp^k.
$
Thus, the reduced PPPA \eqref{eq:forward-main} can be viewed either as a
relaxed proximal iteration for $A_C$, as in
\cite{BrediesEtAl2022}, or as a forward iteration
for the complementary resolvent $F=J_{A_C^{-1}}$.
The latter viewpoint allows us to directly apply the convergence
results of \cref{sec:forward-rel-semi},
since relative semimonotonicity of $A$ yields
that of $F$;
see \cref{lem:transfer-lemma}.

\paragraph{Range-restricted operator.} As observed by
\cite[Eq.~(2.4) and Proposition~2.5]{BrediesEtAl2022}
and further explored by
\cite[Lemma~2.2 and Definition~2.1]{XueZhang2026},
only graph points of $A$ whose operator values
belong to $\Rr$ are relevant to the
preconditioned resolvent. To make this explicit,
define
\begin{equation}
	A_{\Rr}(x)\coloneqq Ax\cap\Rr,
	\qquad x\in\Hh. \label{eq:A-Rr}
\end{equation}
Since $Cz-P\bar{x}\in\Rr$ for all $z\in\Rr$
and $\bar{x}\in\Hh$, we have
$
(P+A_{\Rr})^{-1}C
=(P+A)^{-1}C.
$
Moreover, $0\in\Rr$ and $\ran C=\Rr$
imply that
$
\zer A_{\Rr}=\zer A$ and $
(A_{\Rr})_C=A_C.
$
In particular, the reduced forward map admits
the equivalent representation
$
F=J_{((A_{\Rr})_C)^{-1}}.
$
Thus, replacing $A$ with $A_{\Rr}$ changes
neither the PPPA nor its reduced forward
representation. As will be seen below, the
semimonotonicity and demiclosedness assumptions can be imposed on
$A_{\Rr}$ instead of on $A$. Similarly, sufficient
error bounds involving the original operator
can be formulated using $A_{\Rr}$.

\subsection{Basic convergence and projected solution identification}
\label{subsec:relative-basic}

The following lemma transfers relative semimonotonicity of $A_\Rr$ given by \eqref{eq:A-Rr}
to $F$.

\begin{lemma}[Transfer of relative semimonotonicity]
	\label{lem:transfer-lemma}
Let $\varnothing\ne S^\star\subseteq\zer A$, and suppose
that $A_{\Rr}$ is relatively $(U,V)$-semimonotone
with respect to $S^\star$, that is,
\[
\langle x-x^\star,a\rangle
\ge
\langle x-x^\star,U(x-x^\star)\rangle
+\langle a,Va\rangle \quad \forall x^\star\in S^\star, ~x\in\Hh, ~a\in Ax\cap\Rr,
\]
and let $F$ be given by \eqref{eq:reduced-forward-map}. Suppose that
	$U\succeq\tau P$ for some $\tau\in \R$. Let $q\in F(z)$ and $\bar{x}\in(P+A)^{-1}Cz$ satisfy
	$q=z-C^*\bar{x}$. Then, for every $x^\star\in S^\star$,
	with $d\coloneqq z-C^*x^\star$,
	\begin{equation}
		(1+2\tau)\ip{q}{d}
		\ge
		\tau\norm{d}^2
		+\ip{q}{(H+\tau\Id_{\Rr})q}
		+\norm{(U-\tau P)^{1/2}(\bar{x}-x^\star)}^2,
		\label{eq:transfer-master}
	\end{equation}
	where $H\coloneqq \Id_{\Rr}+C^*VC$.
	Consequently, if $\tau >-1/2$, $F$ is relatively
	$\bigl(\frac{\tau}{1+2\tau}\Id_{\Rr},
	\frac{H+\tau\Id_{\Rr}}{1+2\tau}\bigr)$-semimonotone
	with respect to $C^*S^\star$, where $C^* S^\star \subseteq \zer F$. 
\end{lemma}

\begin{proof}
	Since $Cq\in A\bar{x}\cap \Rr$ and $C^*(\bar{x}-x^\star)=d-q$,
	relative semimonotonicity gives
	\[
	\ip{q}{d}
	=\norm{q}^2+\ip{q}{C^*(\bar{x}-x^\star)}
	\ge
	\ip{q}{Hq}
	+\ip{\bar{x}-x^\star}{U(\bar{x}-x^\star)}.
	\]
	Writing $U=\tau P+(U-\tau P)$ and using
	$\ip{\bar{x}-x^\star}{P(\bar{x}-x^\star)}=\norm{d-q}^2$,
	we obtain
	\[
	\ip{q}{d}
	\ge
	\ip{q}{Hq}
	+\tau\norm{d-q}^2
	+\norm{(U-\tau P)^{1/2}(\bar{x}-x^\star)}^2,
	\]
	from which we get
	\eqref{eq:transfer-master}. Since $U-\tau P\succeq0$
	and $1+2\tau>0$, dropping the last term and dividing
	by $1+2\tau$ proves the semimonotonicity claim on $F$. Finally, $C^*S^\star \subseteq \zer F$ follows from \eqref{eq:zeros-F}. 
\end{proof}

In the simple case $U\succeq0$, taking $\tau=0$ in
\cref{lem:transfer-lemma} shows that $F$ is relatively $(0,H)$-semimonotone
with respect to $C^*S^\star$.
We can then apply \cref{prop:forward-convergence} directly to
\eqref{eq:forward-main}.

\begin{corollary}[Convergence of the reduced PPPA]
	\label{cor:pppa-forward-direct}
	Let $\varnothing\ne S^\star\subseteq\zer A$, and suppose that
	$A_\Rr$ is relatively $(U,V)$-semimonotone with respect to $S^\star$,
	where $U\succeq0$ and $H$ is strongly monotone.
	Let $\{(x^k,\bar x^k)\}$ be generated by \eqref{eq:pppa}, and let $ (z^k,q^k)$ be given by \eqref{eq:zq}. Suppose that
	\begin{equation}
		\lambda_k\in (0,2\beta_k] \text{ where } \beta_k>0 \text{ satisfies }
		\ip{q^k}{Hq^k}\ge\beta_k\norm{q^k}^2.
		\label{eq:lambda-wm}
	\end{equation}
	Set $\delta_k\coloneqq \lambda_k(2\beta_k-\lambda_k)$.
	Then the following hold.
	\begin{enumerate}[(i)]
		\item $\{z^k\}$ is Fej\'er monotone with respect to
		$C^*S^\star$; that is, $\{x^k\}$ is Fej\'er monotone with respect to
		$S^\star$ in the $P$-seminorm.  Moreover, $\{q^k\}$,
		$\{\Pi_{\Rr}x^k\}$, and $\{\Pi_{\Rr}\bar x^k\}$ are bounded.
		
		\item If $\liminf_{k\to\infty}\delta_k>0$, then
		$\sum_{k=0}^\infty\norm{q^k}^2<+\infty$,
		$q^k\to0$, and $P(x^k-\bar x^k)\to0$.
		
		\item For every $x^\star\in S^\star$ and $N\in\N$,
		$\norm{\frac{\sum_{k=0}^N\lambda_kq^k}
			{\sum_{k=0}^N\lambda_k}}
		\le 2\norm{x^0-x^\star}_P/\sum_{k=0}^N\lambda_k$,
		which is $O(N^{-1})$ if $\liminf_k\lambda_k>0$.
		If $\delta\coloneqq \inf_k\delta_k>0$, then also
		$\norm{(N+1)^{-1}\sum_{k=0}^Nq^k}
		\le\norm{x^0-x^\star}_P/\sqrt{\delta(N+1)}$.
		
		\item For every $x^\star\in S^\star$ and $N$ such that
		$\sum_{k=0}^N\delta_k>0$,
		$
		\min_{0\le k\le N}\norm{q^k}^2
		=\min_{0\le k\le N}\norm{x^k-\bar x^k}_P^2
		\le\frac{\norm{x^0-x^\star}_P^2}
		{\sum_{k=0}^N\delta_k}.
		$
		If $\delta_k\ge\kappa>0$ for all $k$, the right-hand
		side is at most $\norm{x^0-x^\star}_P^2/[\kappa(N+1)]$.
	\end{enumerate}
\end{corollary}

\begin{proof}
	By \cref{lem:transfer-lemma} and \cref{eq:forward-main}, we can directly invoke \cref{prop:forward-convergence}. From \cref{prop:forward-convergence}(i), we obtain the claim on $\{z^k\}$. Using the definition of $z^k$, we get that $\norm{C^*x^{k+1} -C^*x^\star}^2\leq \norm{C^*x^k - C^*x^\star}^2$ for all $k\in \mathbb{N}$, which proves the Fej\'er monotonicity on $\{x^k\}$ in the $P$-seminorm. Moreover, we have from \eqref{eq:transfer-master} with $\tau=0$ that $\ip{q^k}{z^k-C^*x^\star} \geq \ip{q^k}{Hq^k} \geq \eta_H \norm{q^k}^2$, where 
	\begin{equation}
			\eta_H\coloneqq
		\inf_{\norm{q}=1}\ip{q}{Hq}
		=\norm{H^{-1}}^{-1},
		\label{eq:eta-H}
	\end{equation}
	which is positive since $H$ is strongly monotone.  Thus, 
	$\eta_H\norm{q^k}\le\norm{z^k-C^*x^\star}$ and we get the boundedness of $\{q^k\}$.
	With the boundedness of $\{z^k\}$ and $\{q^k\}$, while noting \eqref{eq:graph-rel} and \eqref{eq:proj-via-C}, we have 
	\begin{equation}
		\Pi_{\Rr}x^k=(P|_{\Rr})^{-1/2}z^k \quad \text{and} \quad 
		\Pi_{\Rr}\bar x^k=(P|_{\Rr})^{-1/2}(z^k-q^k),
		\label{eq:projected-ppa-iterates}
	\end{equation}
	we get the rest of the claim of (i). 
	Finally,  parts (ii)--(iv) directly follow from \cref{prop:forward-convergence}(ii)--(iv) by noting  \cref{eq:graph-rel,eq:residual-equality} and
	$\norm{z^0-C^*x^\star}=\norm{x^0-x^\star}_P$.
\end{proof}

\paragraph{Projected identification and convergence.}
The estimates above do not by themselves identify the weak cluster
points of $\{z^k\}$ as zeros of $F$. This requires a demiclosedness
property of $F$.

\begin{proposition}[Projected identification and convergence]
	\label{prop:projected-identification}
	Suppose the assumptions of \cref{cor:pppa-forward-direct} hold,
	$\liminf_k\delta_k>0$, and $F$ is demiclosed at zero.
	Then the following hold.
	\begin{enumerate}[(i)]
		\item Every weak cluster point of $\{z^k\}$ belongs to
		$C^*(\zer A)$. Consequently, every weak cluster point of
		$\{\Pi_{\Rr}x^k\}$ or $\{\Pi_{\Rr}\bar x^k\}$
		belongs to $\Pi_{\Rr}(\zer A)$.
		
		\item If, in addition,
		\begin{equation}
			C^*S^\star=C^*(\zer A),
			\label{eq:projected-S-equals}
		\end{equation}
		then $z^k\weakto z^\infty\in C^*(\zer A)$,
		and both 	$\{\Pi_{\Rr}x^k\}$ or $\{\Pi_{\Rr}\bar x^k\}$ converge weakly to
		$(P|_{\Rr})^{-1/2}z^\infty\in\Pi_{\Rr}(\zer A)$.
	\end{enumerate}
\end{proposition}

\begin{proof}
	By \cref{prop:forward-convergence}(ii), $q^k\to 0$ and every weak cluster
	point of $\{z^k\}$ belongs to $\zer F=C^*(\zer A)=(P|_{\Rr})^{1/2}\Pi_\Rr (\zer A)$, where the last equality holds by \eqref{eq:C-canonical}. Consequently, by \cref{eq:projected-ppa-iterates}, weak cluster points of $\Pi_{\Rr}x^k $ and
	$\Pi_{\Rr}\bar x^k$ belong to the set  $(P|_{\Rr})^{-1/2}\zer F = \Pi_\Rr (\zer A)$, proving (i). 
	Under \eqref{eq:projected-S-equals}, the weak convergence
	follows again from \cref{prop:forward-convergence}(ii), \eqref{eq:projected-ppa-iterates} and \eqref{eq:C-canonical}.
\end{proof}

Global $(U,V)$-semimonotonicity provides a sufficient
condition for the demiclosedness assumption in
\cref{prop:projected-identification}, together with
a last-iterate residual rate.

\begin{proposition}[Convergence under global $(U,V)$-semimonotonicity]
	\label{prop:global-semimonotonicity}
	Suppose that $\zer A\ne\varnothing$, that
	$A_{\Rr}$ is globally $(U,V)$-semimonotone
	with $U\succeq0$, and that $H$ is strongly
	monotone. Let $\eta_H$ be given by \eqref{eq:eta-H}. Then $F$ is $\eta_H$-cocoercive,
	and hence single-valued, $\eta_H^{-1}$-Lipschitz
	continuous, and demiclosed at zero.
	
	If $0\le\lambda_k\le2\eta_H$ and
	$\sum_{k=0}^\infty
	\lambda_k(2\eta_H-\lambda_k)=+\infty$,
	then $\{z^k\}$ is Fej\'er monotone with respect
	to $C^*(\zer A)$, $q^k\to0$, and
	$z^k\weakto z^\infty\in C^*(\zer A)$.
	Moreover,
	$\Pi_{\Rr}x^k,\Pi_{\Rr}\bar x^k
	\weakto(P|_{\Rr})^{-1/2}z^\infty
	\in\Pi_{\Rr}(\zer A)$.
	If, in addition,
	$\inf_k\lambda_k(2\eta_H-\lambda_k)>0$,
	then $\norm{x^k-\bar x^k}_P^2=o(k^{-1})$.
\end{proposition}

\begin{proof}
	Let $(z,q),(z',q')\in\gra F$, and choose
	$\bar x,\bar y\in\Hh$ such that
	$Cq\in A_{\Rr}\bar x$, $Cq'\in A_{\Rr}\bar y$,
	$C^*\bar x=z-q$, and $C^*\bar y=z'-q'$.
	By global semimonotonicity,
	$
	\ip{q-q'}{z-z'}
	\ge\ip{q-q'}{H(q-q')}
	\ge\eta_H\norm{q-q'}^2.
	$
	Hence, $F$ is $\eta_H$-cocoercive and,
	by the standing full-domain assumption,
	single-valued and $\eta_H^{-1}$-Lipschitz
	continuous. By \cite[Proposition~4.39]{BauschkeCombettes2017},
$N\coloneqq\Id_{\Rr}-2\eta_HF$ is nonexpansive.
Hence, the demiclosedness principle
\cite[Corollary~4.28]{BauschkeCombettes2017}
implies that $F$ is demiclosed at zero. The remaining claims follow from
\cite[Proposition~2.9]{BartzDaoPhan2022}
by expressing \eqref{eq:forward-main} as a
Krasnosel'skii--Mann iteration of $N$, noting that
$\Fix N=\zer F=C^*(\zer A)$, and using
\eqref{eq:projected-ppa-iterates} and
\eqref{eq:residual-equality}.
\end{proof}

\paragraph{Failure of projected identification.}
In contrast to the above result, relative semimonotonicity alone
does not ensure demiclosedness of $F$, even
when $A$ is continuous in finite dimensions.

\begin{example}[Projected limit outside the projected zero set]
	\label{ex:projected-counterexample}
	Let $\Hh=\R^2$, 
	$
	P=\begin{pmatrix}1&0\\0&0\end{pmatrix}$ and 
	$A(u,v)=\left(\frac{u}{1+v^2},0\right).
	$
	Then $A$ is continuous, $\zer A=\{0\}\times\R$,
	and $\ip{(u,v)-(0,t)}{A(u,v)}
	=u^2/(1+v^2)\ge0$ for every $t\in\R$.
	Thus $A$ is relatively $(0,0)$-semimonotone
	with respect to $\zer A$. For $x=(s,t)$, every
	$\bar x=(u,v)\in(P+A)^{-1}Px$ satisfies
	$
	u=s(1+v^2)/(2+v^2),
	$
	where $v\in\R$ is arbitrary. Hence,
	$(P+A)^{-1}P$ has full domain and
	$
		F(s,0)=
		\left\{
		\left(\frac{s}{2+v^2},0\right):
		v\in\R
		\right\}
$ for all $(s,0)\in \Rr$.
	
	Since $H=\Id_{\Rr}$, we may take
	$\beta_k\equiv1$ and $\lambda_k\equiv1$,
	so that  $\delta_k\equiv1$ and \eqref{eq:pppa} reduces to $x^{k+1} \in (P+A)^{-1}Px^k$ and $x^{k+1} = \bar{x}^k$.
	Starting from $x^0=(s_0,0)$ with $s_0>0$,
	the sequence $x^k = (s_k ,k)$ can be generated from
	\eqref{eq:pppa} where
	\begin{equation*}
		s_{k+1}
		=s_k\frac{1+(k+1)^2}{2+(k+1)^2}
		=s_k\left(1-\frac{1}{2+(k+1)^2}\right) = s_0 \prod_{j=0}^{k}(1-a_j),
	\end{equation*}
	with $a_j=(2+(j+1)^2)^{-1}$. Since $\sum_{j=0}^{\infty} a_j<+\infty$, $s_k\to s_\infty>0$.
	Consequently,
	$\Pi_{\Rr}x^k,\Pi_{\Rr}\bar x^k
	\to(s_\infty,0)$, although
	$\Pi_{\Rr}(\zer A)=\{(0,0)\}$.	Meanwhile,
	$
	q^k=
	\left(\frac{s_k}{2+(k+1)^2},0\right)
	\to0,
	$
	while $0\notin F(s_\infty,0)$.
	Thus $F$ is not demiclosed at zero, so that the hypothesis of \cref{prop:projected-identification} does not hold.\qed
\end{example}

Thus, the projected identification and convergence assertions of
\cite[Theorem~2.4(iii),(v)]{EvensEtAl2025} require an additional
assumption. By \cref{prop:projected-identification}, demiclosedness
of $F$ at zero suffices, together with
\eqref{eq:projected-S-equals} for convergence.

\paragraph{Sufficient conditions for demiclosedness of $F$.}
We next characterize demiclosedness of $F$ at zero
and provide a sufficient condition for it.

\begin{lemma}[Demiclosedness and bounded recovery]
	\label{lem:sufficient-demiclosed}
	The following are equivalent.
	\begin{enumerate}[(i)]
		\item $F$ is demiclosed at zero.
		\item $A_C$ is demiclosed at zero.
		\item If $q^j\to0$, $Cq^j\in A_{\Rr}x^j$, and
		$C^*x^j\weakto p$, then $p\in C^*(\zer A)$.
	\end{enumerate}
	Moreover, these conditions hold if $A_{\Rr}$ is
	demiclosed at zero and the bounded-recovery condition (see \cref{def:recovery}(i)) holds.
\end{lemma}

\begin{proof}
	By \cref{prop:parallel-composition}, $q\in F(z)$ if and
	only if $q\in A_C(z-q)$, and $\zer F=\zer A_C$ by
	\eqref{eq:zeros-F}. If $q^j\to0$, then $z^j\weakto z$
	if and only if $z^j-q^j\weakto z$. Hence (i) and (ii)
	are equivalent.
	Since $A_C^{-1}=C^*A^{-1}C$, we have $q^j\in A_Cp^j$
	if and only if $p^j=C^*x^j$ for some $x^j\in\Hh$ with
	$Cq^j\in Ax^j$, equivalently $Cq^j\in A_{\Rr}x^j$,
	because $Cq^j\in\Rr$. As $\zer A_C=C^*(\zer A)$,
	(ii) and (iii) are equivalent.
	
	Finally, suppose that $A_{\Rr}$ is demiclosed at zero
	and the bounded-recovery condition holds. Let the sequences
	in the premise of (iii) be given, and set
	$z^j\coloneqq C^*x^j+q^j$. Then
	$x^j\in(P+A)^{-1}Cz^j$, and $\{z^j\}$ is bounded
	because $z^j\weakto p$. By the bounded-recovery condition,
	$\{x^j\}$ is bounded, so $x^j\weakto x$ along a
	subsequence. Since $Cq^j\in A_{\Rr}x^j$ and
	$Cq^j\to0$, demiclosedness of $A_{\Rr}$ gives
	$x\in\zer A$. Moreover, $C^*x=p$, so (iii) holds.
\end{proof}

Global Lipschitz continuity of $(P+A)^{-1}$
implies the bounded-recovery condition and was
assumed in \cite[Theorem~2.9]{BrediesEtAl2022}
for full-space convergence under monotonicity.
Assuming that $P$ is admissible for $A$,
Xue--Zhang
\cite[Assumption~7 and Proposition~5.1]{XueZhang2026}
weakened this requirement to Lipschitz
continuity of $(P+A)^{-1}|_{\Rr}$.
Both assumptions imply the bounded-recovery
condition, which requires neither Lipschitz
continuity nor single-valuedness.

\subsection{Convergence of the PPPA sequences}
\label{subsec:pppa-convergence}

We now turn to the full PPPA sequences $\{x^k\}$ and
$\{\bar x^k\}$. The following lemma identifies the weak cluster
points of the shadows and gives sufficient conditions for weak
convergence of both sequences.

\begin{lemma}[Weak convergence of the PPPA sequences]
	\label{lem:shadow-weak}
	Let $\{(x^k,\bar x^k)\}$ be generated by
	\eqref{eq:pppa}, with $z^k,q^k$ defined by
	\eqref{eq:zq}. Suppose that $q^k\to0$ and
	$A_{\Rr}$ is demiclosed at zero. Then every
	weak cluster point of $\{\bar x^k\}$ belongs
	to $\zer A$. If, in addition, $z^k\weakto z^\infty$, let
	\begin{equation}
			X\coloneqq
		\{x\in\zer A:C^*x=z^\infty\}.
		\label{eq:X-zerA}
	\end{equation}
	Then the following hold.
	\begin{enumerate}[(i)]
		\item Every weak cluster point of
		$\{\bar x^k\}$ belongs to $X$.

	\item If $\{\bar x^k\}$ is bounded, then
	$X\ne\varnothing$. If, moreover,
	$X=\{x^\infty\}$, then
	$\bar x^k\weakto x^\infty$.
	In particular, $X$ is a singleton whenever
	$P$ is admissible for $A$.
	
		\item If $\bar x^k\weakto x^\infty$
		and
		\begin{equation}
			0<\underline\lambda\le\lambda_k
			\le\overline\lambda<2
			\qquad(k\in\N),
			\label{eq:full-iterate-relaxation}
		\end{equation}
		then $x^k\weakto x^\infty$.
	\end{enumerate}
\end{lemma}

\begin{proof}
	Let $\bar x^{k_j}\weakto\widehat x$.
	By \eqref{eq:graph-rel},
	$Cq^{k_j}\in A_{\Rr}\bar x^{k_j}$
	and $Cq^{k_j}\to0$.
	Demiclosedness of $A_{\Rr}$ gives
	$\widehat x\in\zer A_{\Rr}=\zer A$.
	If $z^k\weakto z^\infty$, then
	$
	C^*\bar x^{k_j}
	=z^{k_j}-q^{k_j}\weakto z^\infty,
	$
	whence $C^*\widehat x=z^\infty$.
	This proves (i). 
	
	For (ii), boundedness of $\{\bar x^k\}$
	ensures the existence of weak cluster
	points, so (i) gives $X\ne\varnothing$.
	If $X=\{x^\infty\}$, every weakly
	convergent subsequence of $\{\bar x^k\}$
	has limit $x^\infty$. By boundedness,
	it follows that
	$\bar x^k\weakto x^\infty$. Now, suppose that $P$ is admissible for $A$.
	For every $x\in X$, we have
	$x\in\zer A$ and $C^*x=z^\infty$.
	Hence, 
	$
	Cz^\infty=Px\in(P+A)x,
	$
	so $x\in(P+A)^{-1}Cz^\infty$.
	Admissibility implies that $X$ has at
	most one element. Since $X\ne\varnothing$
	under the boundedness assumption,
	the conclusion follows.

Finally, suppose that $\bar x^k\weakto x^\infty$
and \eqref{eq:full-iterate-relaxation} holds.
Fix $h\in\Hh$ and set
$
a_k\coloneqq\ip{x^k-x^\infty}{h}$ and 
$b_k\coloneqq\ip{\bar x^k-x^\infty}{h}.$
Since $b_k\to0$ and
$
\theta\coloneqq
\max\{|1-\underline\lambda|,
|1-\overline\lambda|\}<1,
$
the recursion
$a_{k+1}=(1-\lambda_k)a_k+\lambda_kb_k$ derived from \eqref{eq:pppa}
gives
$
|a_{k+1}|
\le\theta|a_k|+\overline\lambda|b_k|.
$
By \cite[\S~2.2.1, Lemma~3]{Polyak1987},
$a_k\to0$. Since $h\in\Hh$ was arbitrary,
$x^k\weakto x^\infty$, proving (iii).
\end{proof}

We now combine the reduced convergence results with
\cref{lem:shadow-weak} to establish projected
convergence and, under additional conditions,
weak convergence of the full PPPA sequences.

\begin{theorem}[Weak convergence of the PPPA sequences]
	\label{thm:all-sequences}
	Suppose the assumptions of \cref{cor:pppa-forward-direct}
	hold, $\liminf_k\delta_k>0$, and
	\eqref{eq:projected-S-equals} holds. Assume further that
	$A_{\Rr}$ is demiclosed at zero and the
	bounded-recovery condition holds. Then there exists
	$z^\infty\in C^*(\zer A)$ such that
	\begin{equation}
			z^k\weakto z^\infty,\qquad
		\Pi_{\Rr}x^k,\,\Pi_{\Rr}\bar x^k
		\weakto (P|_{\Rr})^{-1/2}z^\infty
		\in\Pi_{\Rr}(\zer A). 
		\label{eq:reduced-and-projected-convergence}
	\end{equation}
	Moreover, $\{\bar x^k\}$ is bounded, and
	the set $X\subseteq \zer A$ defined by \eqref{eq:X-zerA}
	is nonempty and contains every weak cluster
	point of $\{\bar x^k\}$.
	If $X=\{x^\infty\}$, then
	$\bar x^k\weakto x^\infty$.
	This uniqueness condition holds, in particular,
	when $P$ is admissible for $A$.
	If also \eqref{eq:full-iterate-relaxation} holds,
	then $x^k\weakto x^\infty$.
\end{theorem}

\begin{proof}
	By \cref{cor:pppa-forward-direct}, $\{z^k\}$
	is bounded and $q^k\to0$.
	Demiclosedness of $A_{\Rr}$ at zero and the
	bounded-recovery condition imply demiclosedness
	of $F$ at zero by \cref{lem:sufficient-demiclosed}.
	Thus, \cref{prop:projected-identification}(ii)
	yields \eqref{eq:reduced-and-projected-convergence}.
	Since $\{z^k\}$ is bounded, the bounded-recovery
	condition also implies boundedness of
	$\{\bar x^k\}$. The remaining conclusions
	follow from \cref{lem:shadow-weak}.
\end{proof}

\subsection{Further consequences of the monotonicity modulus}
\label{subsec:pppa-modulus}
We now exploit the monotonicity modulus $U$ of $A_{\Rr}$. First,
the nonnegative term in \eqref{eq:transfer-master} yields a
Lipschitz recovery estimate (\cref{def:recovery}(ii)).

\begin{theorem}[Lipschitz recovery from the monotonicity modulus]
	\label{thm:kernel-lift}
	Let $\varnothing\ne S^\star\subseteq\zer A$, and suppose that
	$A_{\Rr}$ is relatively $(U,V)$-semimonotone with respect to
	$S^\star$. Suppose that $U\succeq\tau P$ for some $\tau\in\R$,
	$H+\tau\Id_{\Rr}$ is strongly monotone, and $U$ is strongly
	monotone on $\ker P$, that is,
	$\ip{h}{Uh}\ge\gamma\norm{h}^2$ for all $h\in\ker P$
	and some $\gamma>0$. 
Then the Lipschitz recovery estimate
\eqref{eq:lipschitz-recovery} holds relative to
$S^\star$. Consequently, the bounded-recovery
condition holds and
$\Phi(C^*x^\star)=\{x^\star\}$
for every $x^\star\in S^\star$.
\end{theorem}

\begin{proof}
	Fix $x^\star\in S^\star$, $z\in\Rr$, and
	$\bar x\in(P+A)^{-1}Cz$. Set
	$d\coloneqq z-C^*x^\star$, $q\coloneqq z-C^*\bar x$,
	and $W\coloneqq U-\tau P\succeq0$. Let $\eta>0$ be such that
	$H+\tau\Id_{\Rr}\succeq\eta\Id_{\Rr}$.
	By \eqref{eq:transfer-master} and Cauchy--Schwarz,
	\[
	\norm{W^{1/2}(\bar x-x^\star)}^2+\eta\norm{q}^2
	\le (1+2\tau)\ip{q}{d}-\tau\norm{d}^2
	\le |1+2\tau|\norm{q}\norm{d}
	+\max\{-\tau,0\}\norm{d}^2.
	\]
	Applying Young's inequality to the last expression yields
	$
	\norm{W^{1/2}(\bar x-x^\star)}^2
	+\frac{\eta}{2}\norm{q}^2\le L\norm{d}^2,
	$
	where $L\coloneqq\max\{-\tau,0\}+(1+2\tau)^2/(2\eta)$.
	Hence, $\norm{q}\le c\norm{d}$ and
	$\norm{W^{1/2}(\bar x-x^\star)}\le\sqrt L\norm{d}$,
	where $c\coloneqq\sqrt{2L/\eta}$.

	Write $\bar x-x^\star=r+h$ with $r\in\Rr$ and
	$h\in\ker P$. Then $C^*r = C^*(\bar x- x^\star) = d - q$ and by \eqref{eq:C-canonical}, $(P|_\Rr)^{1/2}r = d - q$. Thus,
	$\norm r\le M(\norm d + \norm q) \le M(1+c) \norm d$, where 	$M\coloneqq\norm{(P|_{\Rr})^{-1/2}}$. On the other hand, $\ip{h}{Wh}=\ip{h}{Uh}\ge\gamma\norm h^2$, and so 
	$\sqrt\gamma\norm h\le\norm{W^{1/2}h}
	\le\sqrt L\norm d+\sqrt{\norm W}\norm r$.
	Combining these estimates proves \eqref{eq:lipschitz-recovery} with
	$K\coloneqq M(1+c)(1+\sqrt{\norm W/\gamma})
	+\sqrt{L/\gamma}$. The rest of the claims follows from \cref{lem:recovery-consequences}.
\end{proof}

If $A_{\Rr}$ is demiclosed at zero, then
\cref{thm:kernel-lift,lem:sufficient-demiclosed}
imply that $F$ is demiclosed at zero.
Moreover, under the assumptions of
\cref{thm:all-sequences} other than bounded-recovery condition,
strong monotonicity of $U$ on $\ker P$ suffices
to obtain $\bar x^k\weakto x^\infty\in S^\star$,
\textit{without requiring admissibility}.

We next apply \cref{prop:forward-convergence}(v) to directly obtain strong and linear convergence.

\begin{theorem}[Strong and linear convergence of the PPPA]
	\label{thm:weak-minty}
	Let $\varnothing\ne S^\star\subseteq\zer A$, and
	suppose that $A_{\Rr}$ is relatively
	$(U,V)$-semimonotone with respect to $S^\star$,
	where $U\succeq\tau P$ for some $\tau>0$ and
	$H+\tau\Id_{\Rr}$ is strongly monotone.
	Suppose that $\lambda_k\in(0,2\beta_k^\tau]$,
	where $\beta_k^\tau>0$ satisfies
	$\ip{q^k}{(H+\tau\Id_{\Rr})q^k}
	\ge(1+2\tau)\beta_k^\tau\norm{q^k}^2$.
	Then $C^*(\zer A)=\{z^\star\}$. Moreover,
	\begin{enumerate}[(i)]
		\item If $\sum_k\lambda_k=+\infty$, then
		$z^k\to z^\star$, $q^k\to0$, and
		$\Pi_{\Rr}x^k,\Pi_{\Rr}\bar x^k
		\to(P|_{\Rr})^{-1/2}z^\star$
		strongly.
		
		\item If $\liminf_k\lambda_k>0$, then
		$\{z^k\}$ converges $Q$-linearly, while
		$\{q^k\}$, $\{\Pi_{\Rr}x^k\}$, and
		$\{\Pi_{\Rr}\bar x^k\}$ converge
		$R$-linearly to their respective limits
		in (i).
		
	\item If, in addition, a local Lipschitz
	recovery estimate holds at some
	$x^\star\in S^\star$, then
	$\zer A=\{x^\star\}$ and
	$\bar x^k\to x^\star$ strongly whenever
	$\sum_k\lambda_k=+\infty$.
	If $\liminf_k\lambda_k>0$, then
	$\{\bar x^k\}$ converges $R$-linearly.
	If also \eqref{eq:full-iterate-relaxation}
	holds, then $x^k\to x^\star$ $R$-linearly.
	\end{enumerate}
\end{theorem}

\begin{proof}
Fix $x^\star\in S^\star$ and set
$z^\star\coloneqq C^*x^\star$.
By \cref{lem:transfer-lemma}, $F$ is relatively
$(a_\tau\Id_{\Rr},V_\tau)$-semimonotone
with respect to $C^*S^\star$, where
$a_\tau\coloneqq\tau/(1+2\tau)>0$ and
$V_\tau\coloneqq
(H+\tau\Id_{\Rr})/(1+2\tau)\succ0$.
Taking $q=0$ in \eqref{eq:transfer-master}
shows that $\zer F=\{z^\star\}$.
Hence, $C^*(\zer A)=\{z^\star\}$ by
\eqref{eq:zeros-F}.

The strong and $Q$-linear convergence
of $\{z^k\}$ in (i) and (ii) follow directly
from \cref{prop:forward-convergence}(v).
Choose $\eta>0$ such that
$V_\tau\succeq\eta\Id_{\Rr}$.
Relative semimonotonicity gives
$\eta\norm{q^k}^2
\le\ip{q^k}{z^k-z^\star}$,
whence $\norm{q^k}
\le\eta^{-1}\norm{z^k-z^\star}$.
Thus, the convergence and rates of
$\{q^k\}$ follow, and
\eqref{eq:projected-ppa-iterates}
gives the corresponding conclusions
for both projected PPPA sequences.

For (iii), \cref{lem:recovery-consequences}
shows that $x^\star$ is the unique zero of
$A$ satisfying $C^*x=z^\star$.
Since $C^*(\zer A)=\{z^\star\}$, we have
$\zer A=\{x^\star\}$.
By (i), whenever $\sum_k\lambda_k=+\infty$,
we have $z^k\to z^\star$, so the local
Lipschitz recovery estimate gives
$\norm{\bar x^k-x^\star}
\le K\norm{z^k-z^\star}$
for all sufficiently large $k$.
Hence, $\bar x^k\to x^\star$ strongly.
If $\liminf_k\lambda_k>0$, then (ii)
and the same bound yield $R$-linear
convergence of $\{\bar x^k\}$.
Finally, under \eqref{eq:full-iterate-relaxation},
the recursion
$x^{k+1}=(1-\lambda_k)x^k+\lambda_k\bar x^k$
implies $R$-linear convergence of $\{x^k\}$,
by an argument analogous to
\cref{lem:shadow-weak}(iii).
\end{proof}

\subsection{Linear convergence under metric subregularity}
\label{subsec:pppa-linear}

We apply \cref{prop:forward-linear} to the reduced PPPA when
$U\succeq-\sigma P$ for some $0\le\sigma<1/2$.
For $x\in\Hh$, let
$\dist_P(x,\zer A)\coloneqq\inf_{x^\star\in\zer A}\norm{x-x^\star}_P$.
By \eqref{eq:zeros-F},
$\dist_P(x,\zer A)=\dist(C^*x,\zer F)$.

\begin{corollary}[Local linear convergence of the PPPA]
	\label{cor:pppa-linear}
	Let $\varnothing\ne S^\star\subseteq\zer A$ satisfy
	\eqref{eq:projected-S-equals}, and suppose that $A_{\Rr}$
	is relatively $(U,V)$-semimonotone with respect to $S^\star$,
	where $U\succeq-\sigma P$ for some $0\le\sigma<1/2$
	and $H-\sigma\Id_{\Rr}$ is strongly monotone.
	Fix $x^\star\in S^\star$ and set
	$z^\star\coloneqq C^*x^\star$.
	Suppose that, for some $R>0$, one of the following holds:
	\begin{enumerate}[(a)]
		\item $F$ is metrically subregular at $z^\star$
		for zero with modulus $\chi>0$, with the error
		bound holding on $B_R(z^\star)$;
		\item $G=A_C$ is metrically subregular at $z^\star$
		for zero with modulus $\widehat\chi>0$, with
		the error bound holding on $B_R(z^\star)$.
	\end{enumerate}
	Let $\{(x^k,\bar x^k)\}$ be generated by \eqref{eq:pppa},
	with $z^k,q^k$ given by \eqref{eq:zq}.
	Set $\eps\coloneqq\sigma/(1-2\sigma)$ and
	$\beta_\sigma\coloneqq
	(\eta_H-\sigma)/(1-2\sigma)>0$.
	Suppose that
	\begin{equation}
		0<\lambda_k\le2\beta_k,\qquad
		\beta_k\ge\beta_\sigma,\qquad
		\ip{q^k}{(H-\sigma\Id_{\Rr})q^k}
		\ge(1-2\sigma)\beta_k\norm{q^k}^2.
		\label{eq:pppa-linear-relaxation}
	\end{equation}
	Under (a) or (b), respectively, set
	\begin{equation}
		\Delta_k\coloneqq
		\begin{cases}
			\lambda_k(2\beta_k-\lambda_k-2\eps\chi^2),
			& \text{under (a)},\\
			\lambda_k[2\beta_k-\lambda_k
			-2\eps(\widehat\chi+|1-\lambda_k|)^2],
			& \text{under (b)}.
		\end{cases}
		\label{eq:pppa-linear-gap}
	\end{equation}
	Suppose that $\Delta_k\ge\kappa$ for all $k$,
	for some $\kappa>0$.
	Then there exist $r\in(0,R)$, $\rho\in[0,1)$,
	and $\Gamma>0$ such that,
	whenever $\norm{x^0-x^\star}_P<r$, the following hold.
	\begin{enumerate}[(i)]
		\item Let $d_k\coloneqq\dist_P(x^k,\zer A)$.
		Then $d_{k+1}\le\rho d_k$, and hence
		$d_k\le d_0\rho^k$.
		Moreover, $\norm{x^k-\bar x^k}_P=\norm{q^k}
		\le m_\eps(\beta_k)d_k$, and
		$\norm{x^{k+1}-x^k}_P\le(1-\rho)\Gamma d_k$.
		
		\item The reduced PPPA sequence $\{z^k\}$ converges strongly
		to some $z^\infty\in C^*(\zer A)$ and, for every $m\ge k$,
		$\norm{x^m-x^k}_P\le\Gamma d_k$,
		$\norm{z^k-z^\infty}\le\Gamma d_k$, and
		$\norm{C^*\bar x^k-z^\infty}
		\le[\Gamma+m_\eps(\beta_k)]d_k$.
		Consequently, both projected PPPA sequences converge
		$R$-linearly to $(P|_{\Rr})^{-1/2}z^\infty$.
		
		\item Suppose, in addition, that
		\eqref{eq:lipschitz-recovery} holds with a constant
		$K>0$ uniformly over $S^\star$.
		Then there exists $x^\infty\in S^\star$,
		the unique zero of $A$ with
		$C^*x^\infty=z^\infty$, such that
		$\norm{\bar x^k-x^\infty}\le K\Gamma d_k$.
		Thus, $\{ \bar x^k\}$ converges strongly and $R$-linearly
		to $x^\infty$ and if \eqref{eq:full-iterate-relaxation}
		holds, then $x^k\to x^\infty$ $R$-linearly.
	\end{enumerate}
	When $\sigma=0$, the bounds in (ii)--(iii) hold with
	$\Gamma_0\coloneqq\min\{2,\Gamma\}$ in place of $\Gamma$,
	and $m_0(\beta_k)=1/\beta_k$.
\end{corollary}

\begin{proof}
By \cref{lem:transfer-lemma} with $\tau=-\sigma$
and \eqref{eq:projected-S-equals}, $F$ is relatively
$(-\eps\Id_{\Rr},V_\sigma)$-semimonotone
with respect to $\zer F=C^*(\zer A)$, where
$V_\sigma\coloneqq(H-\sigma\Id_{\Rr})/(1-2\sigma)
\succeq\beta_\sigma\Id_{\Rr}$.
The standing full-domain assumption ensures that
$F$ has nonempty values everywhere.
By \eqref{eq:forward-main}, the reduced PPPA is
precisely the forward iteration \eqref{eq:forward}
for $F$. Moreover, \eqref{eq:pppa-linear-relaxation}
gives the required coefficients $\beta_k$, while
\eqref{eq:pppa-linear-gap} agrees with \eqref{eq:local-gap}.
Thus, \cref{prop:forward-linear} applies when
$\norm{z^0-z^\star}=\norm{x^0-x^\star}_P<r$
with sufficiently small $r>0$.
Using \cref{eq:residual-equality,eq:projected-ppa-iterates} with \cref{prop:forward-linear}(ii)--(iii)
gives the estimates in (i) and (ii).

For (iii), choose $x^\infty\in S^\star$ with
$C^*x^\infty=z^\infty$, which is possible by
\eqref{eq:projected-S-equals}.
\cref{eq:lipschitz-recovery} gives
$\norm{\bar x^k-x^\infty}
\le K\norm{z^k-z^\infty}\le K\Gamma d_k$.
As in the proof of \cref{thm:weak-minty}(iii),
applying the same estimate at $z=z^\infty$
shows that $x^\infty$ is the unique zero of $A$
with $C^*x^\infty =z^\infty$.
The $R$-linear convergence of $\{x^k\}$ under
\eqref{eq:full-iterate-relaxation} follows
by the same argument as in
\cref{thm:weak-minty}(iii).
Lastly, the improved bounds when $\sigma=0$
follow from \cref{prop:forward-linear}.
\end{proof}

\paragraph{Sufficient conditions for metric subregularity of $F$ and $G$.}
The Lipschitz recovery estimate \eqref{eq:lipschitz-recovery} used
in \cref{cor:pppa-linear}(iii) follows from
\cref{thm:kernel-lift} if $U$ is strongly monotone
on $\ker P$. The following transfers metric
subregularity of $A_{\Rr}$ to $F$ and $G$, using
only a local version of this estimate.

\begin{lemma}[Transfer of metric subregularity]
	\label{lem:pppa-subreg-transfer}
	Let $x^\star\in\zer A$ and $z^\star\coloneqq C^*x^\star\in \zer F$.
	Suppose that $A_{\Rr}$ is metrically subregular at
	$x^\star$ for zero with modulus $\chi_{\Rr}>0$.
	Suppose also that there exist $K,R_L>0$ such that
	\begin{equation}
		\norm{\bar x-x^\star}\le K\norm{z-z^\star}
		\quad\forall z\in B_{R_L}(z^\star),\quad
		\bar x\in(P+A)^{-1}Cz.
		\label{eq:local-subreg-lift}
	\end{equation}
	Then $F$ and $G\coloneqq(F^{-1}-\Id)^{-1}=A_C$
	are metrically subregular at $z^\star$ for zero,
	with moduli $1+c$ and $c$, respectively, where
	$c\coloneqq\chi_{\Rr}\norm P$.
	If the error bound for $A_{\Rr}$ holds globally,
	then both reduced error bounds hold globally
	without assuming \eqref{eq:local-subreg-lift}.
\end{lemma}

\begin{proof}
	By metric subregularity, there exists $r_A>0$ such that
	\begin{equation}
		\dist(x,\zer A)
		\le\chi_{\Rr}\dist(0,A_{\Rr}x)
		\qquad\forall x\in B_{r_A}(x^\star).
		\label{eq:errorbound_A}
	\end{equation}
	Choose $R_0\in(0,R_L)$ with $KR_0<r_A$. Let $z\in B_{R_0}(z^\star)$ and $q\in F(z)$, and choose
	$\bar x\in(P+A)^{-1}Cz$ with
	$C^*\bar x=z-q$ and $Cq\in A_{\Rr}\bar x$.
	Then $\bar x\in B_{r_A}(x^\star)$ by \eqref{eq:local-subreg-lift}, and so by \eqref{eq:errorbound_A}, we get
	\begin{equation}
		\dist(z-q,\zer F)
		\le\norm{C^*}\dist(\bar x,\zer A)
		\le\chi_{\Rr}\norm{C^*}\norm{Cq}
		\le c\norm q.
		\label{eq:pppa-shadow-error-transfer}
	\end{equation}
	Hence, $\dist(z,\zer F)\le(1+c)\norm q$.
	Taking the infimum over $q\in F(z)$ proves
	metric subregularity of $F$.
	
	To establish that of $G$, choose
	$0<s<\min\{R_0/2,cR_0/2\}$.
	If $p\in B_s(z^\star)$ and $q\in G(p)$,
	then $q\in F(p+q)$. To prove metric subregularity of $G$ with constant $c$, it suffices to show that $\dist(p,\zer F)\le c\norm q$.
	When $\norm q<R_0/2$, we have
	$z\coloneqq p+q\in B_{R_0}(z^\star)$. Thus, 
	\eqref{eq:pppa-shadow-error-transfer} gives
	$\dist(p,\zer F)\le c\norm q$.
	On the other hand, if $\norm q \geq R_0/2$, we have
	$\dist(p,\zer F)\le\norm{p-z^\star}
	<s\le cR_0/2\le c\norm q$. Hence, the metric subregularity of $G$ holds. 
	
	If \eqref{eq:errorbound_A} holds on all of $\Hh$,
	the estimate \eqref{eq:pppa-shadow-error-transfer}
	holds for every $(z,q)\in\gra F$ without using
	\eqref{eq:local-subreg-lift}.
	The corresponding bounds for $F$ and $G$
	are therefore global.
\end{proof}

Since $A_{\Rr}x\subseteq Ax$ and $\zer A_{\Rr}=\zer A$,
metric subregularity of $A$ also implies that of $A_{\Rr}$.
By \cref{lem:pppa-subreg-transfer}, metric subregularity of
$A_{\Rr}$ at $x^\star$ supplies either alternative in
\cref{cor:pppa-linear} whenever Lipschitz recovery holds
at $x^\star$. 

\subsection{Finite-dimensional specialization}
\label{sec:evens}

Suppose that $\Hh=\R^n$. If $(P+A)^{-1}P$ is locally bounded at
every point of $\Rr$, then the bounded-recovery condition holds,
since $Cz=P((P|_{\Rr})^{-1/2}z)$ for $z\in\Rr$ and bounded subsets
of $\Rr$ have compact closure. If, in addition, $A_{\Rr}$ is outer
semicontinuous, it is demiclosed at zero, so $F$ is demiclosed at
zero by \cref{lem:sufficient-demiclosed}. Local boundedness holds,
in particular, when $(P+A)^{-1}P$ is single-valued and continuous on
$\Rr$, and outer semicontinuity of $A$ implies that of $A_{\Rr}$.

\begin{corollary}[Finite-dimensional PPPA convergence]
	\label{prop:recover-evens}
	Suppose that $\Hh=\R^n$ and the assumptions of
	\cref{cor:pppa-forward-direct} hold with
	$\liminf_k\delta_k>0$. Suppose further that
	$A_{\Rr}$ is outer semicontinuous and
	$(P+A)^{-1}P$ is locally bounded at every point
	of $\Rr$. Then $F$ is demiclosed at zero,
	$\{\bar x^k\}$ is bounded, and every accumulation
	point of $\{\bar x^k\}$ belongs to $\zer A$.
	Moreover, every accumulation point of
	$\{\Pi_{\Rr}x^k\}$ or $\{\Pi_{\Rr}\bar x^k\}$
	belongs to $\Pi_{\Rr}(\zer A)$.
	
	If $\Pi_{\Rr}S^\star=\Pi_{\Rr}(\zer A)$,
	then both projected PPPA sequences converge in norm
	to the same point in $\Pi_{\Rr}(\zer A)$.
	If, in addition, $P$ is admissible for $A$,
	then $\bar x^k\to x^\infty\in\zer A$.
	Under \eqref{eq:full-iterate-relaxation},
	we also have $x^k\to x^\infty$.
\end{corollary}

\begin{proof}
	By the preceding observations, $A_{\Rr}$ is
	demiclosed at zero and the bounded-recovery condition
	holds. Since $C^*|_{\Rr}$ is an isomorphism,
	$\Pi_{\Rr}S^\star=\Pi_{\Rr}(\zer A)$ is equivalent to
	$C^*S^\star=C^*(\zer A)$. The conclusions follow
	from \cref{lem:sufficient-demiclosed,prop:projected-identification,thm:all-sequences},
	since weak and strong convergence coincide in $\R^n$.
\end{proof}

To compare with \cite{EvensEtAl2025}, let
$E\in\R^{n\times m}$ have orthonormal columns spanning
$\Rr$, and set $X\coloneqq E^\top PE\succ0$.
Identifying $\Rr$ with $\R^m$ through $E$, we have
$C=EX^{1/2}$ and
\begin{equation*}
	\eta_H=1+\lambda_{\min}(X^{1/2}E^\top VEX^{1/2})
	=1+\lambda_{\min}(E^\top VPE).
\end{equation*}
Here, $E^\top VPE$ is similar to the symmetric matrix
$X^{1/2}E^\top VEX^{1/2}$. Hence, strong monotonicity
of $H$ is equivalent to
\cite[Assumption~I.A4]{EvensEtAl2025}.
For $q^k\ne0$, the largest admissible $\beta_k$ in
\eqref{eq:lambda-wm} is
$\ip{q^k}{Hq^k}/\norm{q^k}^2$, which agrees with the
adaptive coefficient in
\cite[Theorem~2.4]{EvensEtAl2025}.

When $U=0$, \cref{cor:pppa-forward-direct} recovers the
residual and best-iterate estimates of
\cite[Theorem~2.4]{EvensEtAl2025}. Under global $(U,V)$-semimonotonicity with
$U\succeq0$, \cref{prop:global-semimonotonicity}
also provides a Hilbert-space counterpart of the
last-iterate squared-residual rate in
\cite[Theorem~2.11]{EvensEtAl2025}.
For the projected assertions in items~(iii) and~(v) of \cite[Theorem~2.4]{EvensEtAl2025}, which also appear as items~(iii) and~(iv) of its
restatement in \cite[Theorem~2.1]{EvensLatafatPatrinos2025},
the argument passes from limit points of the projected shadows
to limit points of the shadows themselves.
This step can fail when the kernel components of the shadows
are unbounded: \cref{ex:projected-counterexample} satisfies the
stated assumptions, yet its projected shadows converge to a
point outside $\Pi_{\Rr}(\zer A)$.
\Cref{prop:recover-evens} shows that local boundedness of
$(P+A)^{-1}P$ suffices. Local boundedness holds, in particular,
under the continuity assumption that both theorems impose for
convergence of the full shadows.\footnote{The first arXiv versions of
	\cite{EvensEtAl2025,EvensLatafatPatrinos2025}
	imposed continuity assumptions on the relevant
	resolvents, under which these conclusions hold.}
For the Douglas--Rachford and Chambolle--Pock applications of
\cite{EvensEtAl2025,EvensLatafatPatrinos2025}, the required recovery
condition follows from their full-graph semimonotonicity
assumptions, so the projected identification and convergence
conclusions of these applications remain valid (\cref{subsec:cp}).

The linear convergence result of \cite[Theorem~2.13]{EvensEtAl2025}
assumes metric subregularity of $A$ at every zero, convexity of
$\zer A$, and continuity of $(P+A)^{-1}P$, and gives an eventual
rate from arbitrary initial points. For $\sigma=0$,
\cref{cor:pppa-linear} instead gives local reduced and projected
rates without convexity or resolvent continuity. When its
subregularity hypothesis is obtained from that of $A_{\Rr}$ through
\cref{lem:pppa-subreg-transfer}, which uses a local Lipschitz
recovery estimate, alternative~(a) yields a squared contraction
factor of the same form as in \cite[Theorem~2.13]{EvensEtAl2025},
while alternative~(b) can yield a sharper one; see
\cref{prop:forward-linear}.

\section{Primal--dual splitting}
\label{sec:cp}

We now turn to Chambolle--Pock and
Douglas--Rachford splitting. Let $\mathcal X$ and $\mathcal Y$ be real Hilbert spaces,
$A:\mathcal X\rightrightarrows\mathcal X$,
$B:\mathcal Y\rightrightarrows\mathcal Y$, and
$0\ne L\in\mathcal B(\mathcal X,\mathcal Y)$. On
$\mathcal H\coloneqq\mathcal X\oplus\mathcal Y$, define
\begin{equation}
	T_{\rm PD}(u,v)\coloneqq
	\begin{bmatrix}Au+L^*v\\B^{-1}v-Lu\end{bmatrix}
	\quad \text{and} \quad 
	P_{\gamma,\tau}\coloneqq
	\begin{bmatrix}
		\gamma^{-1}\Id_{\mathcal X}&-L^*\\
		-L&\tau^{-1}\Id_{\mathcal Y}
	\end{bmatrix},
	\label{eq:A_P-CP}
\end{equation}
with $\gamma,\tau>0.$ Thus, $\zer T_{\rm PD}$ consists of pairs $(u,v)$ satisfying
$-L^*v\in Au$ and $Lu\in B^{-1}v$. In particular, $u$ solves the composite inclusion
$
0\in Au+L^*B(Lu).
$
Moreover, $P_{\gamma,\tau}\succeq0$ if and only if
$\gamma\tau\norm{L}^2\le1$, and it is strongly monotone
when the inequality is strict. At the critical stepsize $\gamma\tau\norm{L}^2=1$,
the range of $P_{\gamma,\tau}$ need not be closed
in Hilbert spaces. Indeed,
\[
P_{\gamma,\tau}
=Q
\begin{bmatrix}
	\gamma^{-1}\Id_{\mathcal X}&0\\
	0&\gamma S
\end{bmatrix}
Q^*, 
\quad\text{where} \quad 
Q\coloneqq
\begin{bmatrix}
	\Id&0\\-\gamma L&\Id
\end{bmatrix}
\]
and $
S\coloneqq\norm L^2\Id_{\mathcal Y}-LL^*$. 
Since $Q$ is boundedly invertible,
$P_{\gamma,\tau}$ has closed range if and only if
$S$ has closed range.

Write $x^k=(u^k,v^k)$ and
$\bar x^k=(\bar u^k,\bar v^k)$ for the PPPA iterates and
shadows, respectively. Then
$\bar x^k\in(P_{\gamma,\tau}+T_{\rm PD})^{-1}
P_{\gamma,\tau}x^k$ is equivalent to
\begin{equation}
	\bar u^k\in J_{\gamma A}(u^k-\gamma L^*v^k) \quad \text{and} \quad 
	\bar v^k\in J_{\tau B^{-1}}
	\bigl(v^k+\tau L(2\bar u^k-u^k)\bigr).
	\label{eq:cp-pppa-updates}
\end{equation}

Together with $x^{k+1}=x^k+\lambda_k(\bar x^k-x^k)$,
this gives the relaxed Chambolle--Pock iteration
\cite{ChambollePock2011}. Full-domain component
resolvents ensure
$\dom((P_{\gamma,\tau}+T_{\rm PD})^{-1}P_{\gamma,\tau})
=\mathcal H$. Henceforth, we suppose $P\coloneqq P_{\gamma,\tau}\succeq0$
has closed range, and let $\Rr\coloneqq\ran P$ and
$P=CC^*$ be as in \eqref{eq:C-canonical} on $\mathcal H$.
We write $F_{\rm PD}$ for the corresponding reduced forward map.

\subsection{Bounded recovery from the primal--dual structure}

We now exploit the primal--dual structure to obtain a Lipschitz
recovery estimate for $T_{\rm PD}$. Unlike \cref{thm:kernel-lift},
this requires no strong monotonicity of $U$ on $\ker P$, but it
assumes relative semimonotonicity on the full graph of
$T_{\rm PD}$, not only on $(T_{\rm PD})_{\Rr}$.

\begin{theorem}[Primal--dual bounded recovery]
	\label{thm:pd-lift}
	Let $\varnothing\ne S^\star\subseteq\zer T_{\rm PD}$
	and suppose that $T_{\rm PD}$ is relatively
	$(U,V)$-semimonotone with respect to $S^\star$
	(on its full graph), where
	$U\succeq-\sigma P$ for some $\sigma\ge0$, where $P=P_{\gamma,\tau}$ is given by \eqref{eq:A_P-CP}, and
	$V=V_1\oplus V_2$ where $V_1,V_2$ are self-adjoint blocks.
	Suppose that $H\coloneqq\Id_{\Rr}+C^*VC$ satisfies
	\begin{equation}
		H-\sigma\Id_{\Rr}\succeq\kappa\Id_{\Rr}
		\qquad\text{for some }\kappa>0.
		\label{eq:cp-coercive}
	\end{equation}
	Then the Lipschitz recovery estimate
	\eqref{eq:lipschitz-recovery} holds relative to
	$S^\star$; that is, there exists $K_{\rm PD}>0$
	such that, for every $x^\star\in S^\star$,
	$x\in\mathcal H$, and
	$\bar x\in(P+T_{\rm PD})^{-1}Px$,
	\begin{equation}
		\norm{\bar x-x^\star}\le K_{\rm PD}
		\norm{C^*(x-x^\star)}.
		\label{eq:pd-bounded-recovery}
	\end{equation}
	Consequently, the bounded-recovery condition holds
	and $(P+T_{\rm PD})^{-1}C(C^*x^\star)=\{x^\star\}$
	for every $x^\star\in S^\star$.
\end{theorem}

\begin{proof}
	Fix $x^\star=(u^\star,v^\star)\in S^\star$,
	$x=(u,v)\in\mathcal H$, and
	$\bar x=(\bar u,\bar v)\in(P+T_{\rm PD})^{-1}Px$.
	Since $x^\star\in\zer T_{\rm PD}$,
	$
	-L^*v^\star\in Au^\star,
	$ and $Lu^\star\in B^{-1}v^\star.
	$
	By \eqref{eq:cp-pppa-updates},
	$\bar u\in J_{\gamma A}(s)$ and
	$\bar v\in J_{\tau B^{-1}}(t)$, where
	$
	s\coloneqq u-\gamma L^*v,$ and 
	$t\coloneqq v+\tau L(2\bar u-u).
	$
	Set $s^\star\coloneqq u^\star-\gamma L^*v^\star$
	and $t^\star\coloneqq v^\star+\tau Lu^\star$.
	
We first derive bounds for
$\norm{\bar u-u^\star}$ and $\norm{\bar v-v^\star}$
in terms of $\norm{s-s^\star}$ and
$\norm{t-t^\star}$, respectively.
	Since $P+PVP=CHC^*$,
	\eqref{eq:cp-coercive} gives
	\begin{equation}
		P+PVP\succeq(\sigma+\kappa)P.
		\label{eq:pd-P-positive}
	\end{equation}
By a direct calculation using the definition of $P$, diagonal blocks of $P+PVP$ are
	$\gamma^{-1}Q_1$ and $\tau^{-1}Q_2$, where
	$
	Q_1\coloneqq\Id_{\mathcal X}
	+\gamma^{-1}V_1+\gamma L^*V_2L$ and 
	$Q_2\coloneqq\Id_{\mathcal Y}
	+\tau^{-1}V_2+\tau LV_1L^*.
	$
	By \eqref{eq:pd-P-positive},
	\begin{equation}
		\begin{aligned}
			\ip{h}{(Q_1-\sigma\Id_{\mathcal X})h}
			&\ge\kappa\norm h^2 \quad \text{and} \quad 	\ip{g}{(Q_2-\sigma\Id_{\mathcal Y})g}
			\ge\kappa\norm g^2,
		\end{aligned}
		\label{eq:pd-block-positive}
	\end{equation}
	for every $h\in\mathcal X$ and
	$g\in\mathcal Y$. Meanwhile, since $\gamma^{-1}(s-\bar u)\in A\bar u$
	and $Lu^\star\in B^{-1}v^\star$, we have
	\[
	\begin{aligned}
		T_{\rm PD}(\bar u,v^\star)
		&=
		\bigl(A\bar u+L^*v^\star,\,
		B^{-1}v^\star-L\bar u\bigr)\\
		&\ni
		\bigl(\gamma^{-1}(s-\bar u)+L^*v^\star,\,
		Lu^\star-L\bar u\bigr)\\
		&=
		\bigl(\gamma^{-1}( s-s^\star-h),-Lh\bigr).
	\end{aligned}
	\]
	where $h\coloneqq\bar u-u^\star$. Applying relative semimonotonicity to this graph
	point and $x^\star$, and using
	$U\succeq-\sigma P$, the definition of $P$ and $V=V_1\oplus V_2$,
	we obtain
	\[
	\begin{aligned}
		\gamma^{-1}\ip{h}{ s-s^\star-h}
		&\ge -\sigma\gamma^{-1}\norm h^2+\gamma^{-2}
		\ip{ s-s^\star-h}{V_1( s-s^\star-h)}
		+\ip{Lh}{V_2Lh}.
	\end{aligned}
	\]
	Multiplying by $\gamma$, expanding the inner products, and rearranging gives
	\[
	\begin{aligned}
		\ip{h}{(Q_1-\sigma\Id_{\mathcal X})h}
		&\le\ip{h}{ s-s^\star}
		+2\gamma^{-1}\ip{h}{V_1( s-s^\star)}-\gamma^{-1}\ip{ s-s^\star}{V_1( s-s^\star)}.
	\end{aligned}
	\]
	Together with \eqref{eq:pd-block-positive},
	this yields
	\begin{equation}
			\begin{aligned}
			\kappa\norm h^2
			&\le\ip{h}{s-s^\star}
			+2\gamma^{-1}\ip{h}{V_1(s-s^\star)}
			-\gamma^{-1}\ip{s-s^\star}{V_1(s-s^\star)}.
		\end{aligned}
		\label{eq:h-primaldual}
	\end{equation}
	
	Similarly,
	since $\tau^{-1}(t-\bar v)\in B^{-1}\bar v$
	and $-L^*v^\star\in Au^\star$, we have
	$
	T_{\rm PD}(u^\star,\bar v)\ni
	\bigl(L^*g,\tau^{-1}(t-t^\star-g)\bigr).
	$ where $g\coloneqq\bar v-v^\star$.
	The same argument gives
	\begin{equation}
			\begin{aligned}
			\kappa\norm g^2
			&\le\ip{g}{t-t^\star}
			+2\tau^{-1}\ip{g}{V_2(t-t^\star)}-\tau^{-1}
			\ip{t-t^\star}{V_2(t-t^\star)}.
		\end{aligned}
		\label{eq:g-primaldual}
	\end{equation}
	Applying Cauchy--Schwarz to \eqref{eq:h-primaldual} and \eqref{eq:g-primaldual}, we get
	\[
	\begin{aligned}
		\kappa\norm h^2
		&\le(1+2a_1)\norm h\norm{s-s^\star}
		+a_1\norm{s-s^\star}^2,\\
		\kappa\norm g^2
		&\le(1+2a_2)\norm g\norm{t-t^\star}
		+a_2\norm{t-t^\star}^2.
	\end{aligned}
	\]
	where $a_1\coloneqq\gamma^{-1}\norm{V_1}$
	and $a_2\coloneqq\tau^{-1}\norm{V_2}$.
	Solving these quadratic inequalities and plugging in $h=\bar{u}-u^\star$, $g = \bar v - v^\star$ gives
	$K_1,K_2>0$, depending only on
	$\gamma,\tau,\kappa,\norm{V_1}$, and
	$\norm{V_2}$, such that
	\begin{equation}
		\norm{\bar u-u^\star}\le K_1\norm{s-s^\star}
		\quad \text{and} \quad 
		\norm{\bar v-v^\star}\le K_2\norm{t-t^\star}.
		\label{eq:pd-component-bounds}
	\end{equation}

On the other hand, by the definitions of $P$, $s$,
and $t$,
$
P(x-x^\star)=(
	\gamma^{-1}(s-s^\star), ~
	\tau^{-1}(t-t^\star)-2L(\bar u-u^\star)
).
$
Combining this identity with
\eqref{eq:pd-component-bounds}, we obtain
\[
\begin{aligned}
	\norm{\bar x-x^\star}
	&\le\norm{\bar u-u^\star}
	+\norm{\bar v-v^\star}\\
	&\le K_1\norm{s-s^\star}
	+K_2\norm{t-t^\star}\\
	&\le\gamma K_1\norm{P(x-x^\star)}
	+\tau K_2\bigl(
	\norm{P(x-x^\star)}
	+2\norm L\norm{\bar u-u^\star}\bigr)\\
	&\le\widehat K_{\rm PD}\norm{P(x-x^\star)},
\end{aligned}
\]
where
$\widehat K_{\rm PD}\coloneqq
\gamma K_1+\tau K_2(1+2\gamma\norm L K_1).
$
Since $P=CC^*\succeq0$ and
$P^2\preceq\norm P\,P$, we have
$
\norm{P(x-x^\star)}
\le\sqrt{\norm P}\norm{C^*(x-x^\star)}.
$
Thus, \eqref{eq:pd-bounded-recovery} follows with
$K_{\rm PD}\coloneqq\sqrt{\norm P}\,\widehat K_{\rm PD}$. The rest of the claim follows from \cref{lem:recovery-consequences}.
\end{proof}
%
%

\subsection{Consequences for Chambolle--Pock and Douglas--Rachford}\label{subsec:cp}

Although the projected assertions of the general result
\cite[Theorem~2.1]{EvensLatafatPatrinos2025} require an additional
hypothesis, such as demiclosedness of the reduced map at zero
(\cref{ex:projected-counterexample,prop:projected-identification}),
the Chambolle--Pock conclusions of that work hold under its
assumptions. In \cite[Assumption~II]{EvensLatafatPatrinos2025}, the
primal--dual weak Minty condition is imposed on the full graph of
$T_{\rm PD}$ with a block-diagonal $V$. Their Stepsize Rule~I ensures
$P\succeq0$ and $1+\bar\eta>0$, where $\bar\eta$ is defined in their
equation~(5). Since $\eta_H=1+\bar\eta$ by \cref{sec:evens},
\eqref{eq:cp-coercive} holds with $\sigma=0$. Thus,
\cref{thm:pd-lift} gives Lipschitz recovery, including at the
critical stepsize. In finite dimensions, outer semicontinuity of $A$
and $B$ implies demiclosedness of $(T_{\rm PD})_{\Rr}$ at zero, and
hence of $F_{\rm PD}$ by \cref{lem:recovery-consequences,lem:sufficient-demiclosed}.
Consequently,
\cref{cor:pppa-forward-direct,prop:projected-identification}
establish the boundedness and projected identification in
\cite[Theorem~3.1(iii)]{EvensLatafatPatrinos2025}. Under their condition
$\Pi_{\Rr}S^\star=\Pi_{\Rr}(\zer T_{\rm PD})$, equivalently
$C^*S^\star=C^*(\zer T_{\rm PD})$,
\cref{thm:all-sequences,thm:pd-lift} also give convergence in norm
of the projected sequences and full shadows, without the continuity
of the component resolvents assumed in
\cite[Theorem~3.1(iv)]{EvensLatafatPatrinos2025}. The full PPPA
iterates converge as well under \eqref{eq:full-iterate-relaxation}.
Douglas--Rachford splitting is the special case $L=\Id$ and
$\tau=1/\gamma$, as observed in \cite{EvensLatafatPatrinos2025}.
Thus, the conclusions of \cite[Theorem~3.3]{EvensEtAl2025} also
hold, with item~(v) requiring no continuity of the resolvents.

\section{Future directions}

The semimonotonicity calculus for Chambolle--Pock
splitting in \cite{EvensLatafatPatrinos2025}
combines the monotonicity and comonotonicity moduli
of the component operators to obtain a
comonotonicity modulus $V$ for $T_{\rm PD}$,
with $U=0$. Our framework allows a nonzero
monotonicity modulus $U$.
Developing a primal--dual semimonotonicity
calculus that retains and exploits this modulus
may yield stronger convergence guarantees. More broadly, the forward viewpoint developed
in this paper suggests revisiting other splitting
algorithms, such as extragradient
\cite{Korpelevich1976},
forward--backward--forward \cite{Tseng2000},
forward--reflected--backward \cite{MalitskyTam2020},
and more general distributed and primal--dual
splitting frameworks
\cite{DaoTamTruong2026,DaoPhanTamTruong2025}.
Our semimonotone splitting framework
\cite{AlcantaraDaoTakeda2026} also admits a
direct forward interpretation, but its analysis
relies on classical cocoercivity estimates.
An interesting question is whether these algorithms,
possibly through suitable reformulations,
admit forward operators satisfying the relative
semimonotonicity conditions of \cref{sec:forward-rel-semi},
leading to more general convergence guarantees.

\section*{AI-assisted research and writing disclosure}
The research was developed by the author, with AI assistance in the following aspects: ChatGPT suggested the concrete counterexample in \cref{ex:projected-counterexample} and assisted with the literature review; ChatGPT and Claude assisted in refining the manuscript's exposition. The author assumes full responsibility for the manuscript.

\section*{Acknowledgments}
The author thanks Emanuele Naldi for helpful discussions and for observing that, in the monotone setting, the forward scheme of \cite{AlcantaraDaoTakeda2026} admits a complementary-resolvent interpretation through the reduced degenerate PPPA framework of \cite{BrediesEtAl2022}. This observation helped motivate the reverse perspective developed here: analyzing degenerate PPPA as a forward method beyond monotonicity.
The author also thanks Brecht Evens for helpful correspondence.

\setlength{\bibsep}{1pt} 	
\bibliographystyle{abbrvnat}
\interlinepenalty=10000
\bibliography{pppa-forward}

\end{document}